\documentclass{jams-l}
\issueinfo{20}{1}{January}{2007}
\pagespan{99}{148}
\dateposted{September 5, 2006}
\PII{S 0894-0347(06)00546-7}
\copyrightinfo{2006}{American Mathematical Society}
\usepackage{url}

\theoremstyle{plain}
\newtheorem{theorem}{Theorem}

\newtheorem{lemma}[theorem]{Lemma}
\newtheorem{proposition}[theorem]{Proposition}
\newtheorem{conjecture}[theorem]{Conjecture}

\numberwithin{theorem}{section}

\theoremstyle{definition}
\newtheorem{definition}[theorem]{Definition}

\numberwithin{equation}{section}

\newcommand{\R}{{\mathbb R}}
\newcommand{\Z}{{\mathbb Z}}
\newcommand{\C}{{\mathbb C}}
\newcommand{\Q}{{\mathbb Q}}

\newcommand{\F}{{\mathbb F}}
\newcommand{\Oc}{{\mathbb O}}
\newcommand{\Ha}{{\mathbb H}}
\newcommand{\Proj}{{\mathbb P}}
\newcommand{\Pic}{\mathop{\textup{Pic}}}
\newcommand{\vol}{\mathop{\textup{vol}}}

\newcommand{\sgn}{\mathop{\textup{sgn}}}

\newcommand{\ev}{\mathop{\textup{ev}}\nolimits}
\newcommand{\tr}{\mathop{\textup{tr}}\nolimits}
\newcommand{\ft}{\widetilde{f}}
\newcommand{\CodeS}{\mathcal{C}}
\newcommand{\CodeT}{\mathcal{C}'}
\newcommand{\co}{\colon}
\begin{document}
\title{Universally optimal distribution\linebreak[1] of points on spheres}
\author{Henry Cohn}
\address{Microsoft Research,
One Microsoft Way,
Redmond, Washington 98052-6399}
\email{cohn@microsoft.com}
\author{Abhinav Kumar}
\address{Department of Mathematics,
Harvard University,
Cambridge, Massachusetts 02138}
\email{abhinav@math.harvard.edu}
\curraddr{Microsoft Research, One Microsoft Way, Redmond, Washington
98052-6399} 
\email{abhinavk@microsoft.com}
\thanks{The second author was supported by a summer internship
in the Theory Group at Microsoft Research and a Putnam Fellowship
at Harvard University.}
\date{November 1, 2004}


\maketitle
\tableofcontents

\section{Introduction} \label{sec:intro}

There are many beautiful and important configurations of points on
the surface of a sphere, ranging from familiar examples such as
the vertices of an icosahedron to more subtle cases such as the
minimal vectors of the $E_8$ or Leech lattices.  Yudin, Kolushov,
and Andreev have shown in \cite{Y,KY1,KY2,A1,A2} that certain
configurations, including the ones mentioned above, minimize
specific sorts of energy.  In this paper, we prove that in fact
they minimize a far broader class of energies than was previously
known;  we call this property \textit{universal optimality\/} (the
precise definition is given below).  Our proofs apply to a number
of other configurations, enumerated in Table~\ref{tab:sharp}
below, as well as to energy minimization in other compact
two-point homogeneous spaces or Euclidean space.

The setting of this paper is potential energy minimization. If one
constrains some identically charged particles to move on the
surface of a sphere, how will they arrange themselves?  They will
eventually spread out so as to locally minimize their electrical
potential energy, assuming their kinetic energy dissipates
according to a force such as viscosity. We will ignore questions
of dynamics and simply focus on the relative positions of the
particles once they reach a minimum for potential energy.
Furthermore, we will focus on global minima, although particles
can indeed become trapped in a local minimum.

This question leads to a fundamental problem in extremal geometry,
in which we abstract the key points by working in an arbitrary
dimension and replacing the electrical potential with an arbitrary
radial potential function.  Let $f \co (0,4] \to [0,\infty)$ be
any decreasing, continuous function.  Given a finite subset
$\CodeS$ of the unit sphere $S^{n-1} \subset \R^n$, define the
\textit{potential energy\/} (more precisely, the
\textit{$f$-potential energy}\/) of $\CodeS$ to be
$$
\sum_{x,y \in \CodeS, x \ne y} f\big(|x-y|^2\big).
$$
Because each pair $x,y$ is counted in both orders, our potential
energy is twice that from physics.  We normalize it this way
because it makes the formulas prettier.  It is also important to
note that we view potential functions as functions of squared
distance.

For a fixed number $N>1$, how can one choose $\CodeS$ with
$|\CodeS|=N$ so as to minimize the potential energy?  (To ensure
that a minimum exists, when $\lim_{r \to 0+}f(r)$ is finite, we
must allow degenerations in which points in the configuration
coincide.) In most cases the answer is some complicated
configuration, but it is remarkable how many important structures
minimize potential energy.

In $\R^3$ the electrical potential corresponds to
$f(r)=1/r^{1/2}$, and a natural generalization to $\R^n$ with
$n\ge3$ is $f(r) = 1/r^{n/2-1}$.  That is an especially important
power law, because $x \mapsto 1/|x|^{n-2}$ is a harmonic function
on $\R^n\setminus\{0\}$. However, any power law is of interest, as
are more general potentials.

Many researchers have studied potential energy minima. One of the
most beautiful rigorous approaches is due to Yudin in \cite{Y}. He
studied the harmonic power law and used spherical harmonics to
give sharp lower bounds for potential energy for $n+1$ or $2n$
points on $S^{n-1}$ (the optima are the regular simplex and cross
polytope).  This work was extended by Kolushov and Yudin in
\cite{KY1}. For $(n,N)=(8,240)$, they showed that the minimal
vectors in the $E_8$ root lattice are the unique global minimum.
Andreev applied their techniques in \cite{A2,A1} to show that the
Leech lattice minimal vectors are the unique minimum for
$(n,N)=(24,196560)$ and the vertices of the regular icosahedron
are the unique minimum for $(n,N)=(3,12)$.  In the present paper,
we extend this approach to much more general potential functions,
as well as numerous other point configurations (listed in
Table~\ref{tab:sharp}).

The paper \cite{KY2} analyzed arbitrary power laws for the
simplex, cross polytope, and $E_8$ root system in a closely
related maximization problem (maximizing potential energy given an
increasing potential function) and raised the question of how
generally the techniques could be applied.  We will deal with
quite general potential functions. Without loss of generality one
can assume that the potential function is nonnegative (only
distances between $0$ and $2$ occur on the sphere, and one can add
a constant to ensure nonnegativity over that range).  The
potential function must be decreasing, which corresponds to a
repulsive force, and it is natural to require convexity as well.
Our results require a strengthening of these conditions, namely
complete monotonicity. Recall that a $C^\infty$ function $f \co I
\to \R$ on an interval $I$ is \textit{completely monotonic\/} if
$(-1)^k f^{(k)}(x) \ge 0$ for all $x \in I$ and all $k \ge 0$ (see
\cite[p.~145]{W}) and \textit{strictly completely monotonic} if
strict inequality always holds in the interior of $I$. Of course,
derivatives at endpoints of intervals denote one-sided
derivatives.  In this paper, endpoints are irrelevant, because we
use half-open intervals closed on the right, and one can show
using the mean value theorem that complete monotonicity on the
interior implies that it holds also at the right endpoint.

If a completely monotonic function fails to be strictly completely
monotonic, then it must be a polynomial: if $f^{(k)}(x)=0$, then
complete monotonicity implies that $f^{(k)}(y)=0$ for all $y>x$,
in which case $f^{(k)}$ is identically zero because it is analytic
(Theorem~3a in \cite[p.~146]{W}).

All inverse power laws $f(r) = 1/r^s$ with $s>0$ are strictly
completely monotonic on $(0,\infty)$, and they are the most
important cases for our purposes.  Other examples include $f(r) =
e^{-cr}$ with $c>0$.

It might seem preferable to use completely monotonic functions of
distance, rather than squared distance, but squared distance is in
fact more natural than it first appears. It simplifies formulas
appearing in later sections of this paper, it fits into the
general framework described in Section~\ref{sec:otherhomog}, and
it is more general than using distance: it is not hard to prove
that if $r \mapsto f(r^2)$ is completely monotonic on a
subinterval $(a,b)$ of $(0,\infty)$, then $f$ is completely
monotonic on $(a^2,b^2)$, but not vice versa.

We will examine certain special configurations of points on
high-dimensional spheres, which were introduced by Levenshtein in
\cite{L2} (in which he proved that they are optimal spherical
codes). We refer to them as the sharp arrangements, because
everything fits into place so perfectly that many bounds turn out
to be sharp. Each is a spherical code with certain parameters: an
$(n,N,t)$ \textit{spherical code\/} is a set of $N$ points on
$S^{n-1}$ such that no two of them have inner product greater than
$t$ (unless they are equal).  In other words, all angles between
them are at least $\cos^{-1} t$.  Each sharp configuration is also
a spherical $M$-design for some $M$; recall that a
\textit{spherical $M$-design\/} is a finite subset of the sphere
such that every polynomial on $\R^n$ of total degree at most $M$
has the same average over the design as over the entire sphere.

\begin{definition} \label{def:sharp}
A finite subset of the unit sphere $S^{n-1}$ is a \textit{sharp
configuration\/} if there are $m$ inner products between distinct
points in it and it is a spherical $(2m-1)$-design.
\end{definition}

Replacing $2m-1$ with $2m+1$ in Definition~\ref{def:sharp} would
be impossible, for the following reason.  If $t_1,\dots,t_m$ are
the inner products that occur between distinct points in a
configuration and $y$ is a point in the configuration, then the
polynomial
$$
x \mapsto \big(1-\langle x,y\rangle\big) \prod_{i=1}^m
\big(\langle x,y \rangle - t_i\big)^2
$$
of degree $2m+1$ vanishes on the entire configuration.  However,
its integral over the sphere does not vanish, because the
polynomial is nonnegative on the sphere and not identically zero.
Thus, sharp configurations are spherical designs of nearly the
greatest possible strength given the number of distances occurring
in them. (Some but not all sharp configurations are actually
$2m$-designs.)

All sharp arrangements that we know of are listed in
Table~\ref{tab:sharp}, together with the $600$-cell (which is not
sharp, but to which our techniques nevertheless apply). The
columns list the dimension $n$ of the ambient Euclidean space, the
number $N$ of points, the largest $M$ such that the code is a
spherical $M$-design, the inner products other than $1$ that occur
between points in the code, and the name of the code, if any.  If
$t$ denotes the largest inner product, then each of these codes is
the unique $(n,N,t)$ spherical code, except for some of those
listed on the last line of the table (the isotropic subspace
codes). See Appendix~\ref{app} for details about uniqueness.

\begin{table}
\caption{The known sharp configurations, together with the $600$-cell.}\label{tab:sharp}
\begin{center}
\begin{tabular}{ccccc}
$n$ & $N$ & $M$ & Inner products & Name\\
\hline $2$ & $N$ & $N-1$ & $\cos (2\pi j/N)$ ($1 \le j \le
N/2$) & $N$-gon\\
$n$ & $N \le n$ & $1$ & $-1/(N-1)$ & simplex\\
$n$ & $n+1$ & $2$ & $-1/n$ & simplex\\
$n$ & $2n$ & $3$ & $-1, 0$ & cross polytope\\
$3$ & $12$ & $5$ & $-1, \pm 1/\sqrt{5}$ & icosahedron\\
$4$ & $120$ & $11$ & $-1,\pm 1/2, 0,(\pm 1 \pm \sqrt{5})/4$ & $600$-cell\\
$8$ & $240$ & $7$ & $-1,\pm 1/2, 0$ & $E_8$ roots\\
$7$ & $56$ & $5$ & $-1, \pm 1/3$ & kissing\\
$6$ & $27$ & $4$ & $-1/2, 1/4$ & kissing/Schl\"afli\\
$5$ & $16$ & $3$ & $-3/5, 1/5$ & kissing\\
$24$ & $196560$ & $11$ & $-1,\pm 1/2, \pm 1/4, 0$ & Leech lattice\\
$23$ & $4600$ & $7$ & $-1,\pm 1/3, 0$ & kissing\\
$22$ & $891$ & $5$ & $-1/2, -1/8, 1/4$ & kissing\\
$23$ & $552$ & $5$ & $-1, \pm 1/5$ & equiangular lines\\
$22$ & $275$ & $4$ & $-1/4, 1/6$ & kissing\\
$21$ & $162$ & $3$ & $-2/7, 1/7$ & kissing\\
$22$ & $100$ & $3$ & $-4/11, 1/11$ & Higman-Sims\\
$q\frac{q^3+1}{q+1}$ & $(q+1)(q^3+1)$ & $3$ & $-1/q, 1/q^2$ &
{\hskip -8pt}
isotropic subspaces\\
& & {\hskip -22pt}(4 if $q=2$){\hskip -22pt}
& & {\hskip -8pt} ($q$ a prime power)\\
\\
\end{tabular}
\end{center}
\end{table}

The first six configurations listed in Table~\ref{tab:sharp} are
the vertices of certain full-dimensional regular polytopes
(specifically, those regular polytopes whose faces are simplices),
together with lower-dimensional simplices. The next seven cases
are derived from the $E_8$ root lattice in $\R^8$ and the Leech
lattice in $\R^{24}$. The $240$-point and $196560$-point
configurations are the minimal (nonzero) vectors in those
lattices.  In sphere packing terms, they are the \textit{kissing
configurations\/}, the points of tangency in the corresponding
sphere packings.  Each arrangement with the label ``kissing'' is
the kissing configuration of the arrangement above it: each
configuration yields a sphere packing in spherical geometry by
centering congruent spherical caps at the points, with radius as
large as possible without making their interiors overlap.  The
points of tangency on a given cap form a spherical code in a space
of one fewer dimension. In general different spheres in a packing
can have different kissing configurations, but that does not occur
here because each configuration's automorphism group acts
transitively on it. (See Chapter~14 of \cite{CS} for the details
of these configurations.) Some of the kissing configurations are
of independent interest. For example, the $27$ points on $S^5$
form the Schl\"afli arrangement, which is connected as follows
with the classical configuration of $27$ lines on a cubic surface.
Let $X$ be a smooth cubic surface in $\C\Proj^3$.  Its divisor
class group $\Pic(X)$ is a free abelian group of rank $7$,
equipped with a bilinear form of signature $(1,6)$ given by taking
intersection numbers (see Proposition~4.8 in \cite[p.~401]{H}).
The class $h$ of a hyperplane section has self-intersection number
$h^2=3$, so its orthogonal complement $h^\perp$ is negative
definite.  Thus, after changing the sign of the metric, $h^\perp
\otimes_\Z \R$ is naturally a Euclidean space. Orthogonally
projecting the divisor classes of the $27$ lines on $X$ from
$\Pic(X) \otimes_\Z \R$ to $h^\perp \otimes_\Z \R$ yields $27$
points on a sphere, which form the Schl\"afli configuration (up to
scaling).  The sharp configurations with $240$, $56$, or $16$
points have similar interpretations via exceptional curves on del
Pezzo surfaces; see Corollary~25.1.1 and Theorem~26.2 in \cite{Ma}
for details.

The sharp configuration of $552$ points in $\R^{23}$ comes from an
equiangular arrangement of $276$ lines in $\R^{23}$ described by
the unique regular two-graph on $276$ vertices; see Chapter~11 of
\cite{GR}. It can be derived from the Leech lattice $\Lambda_{24}$
as follows. Choose any $w \in \Lambda_{24}$ with $|w|^2=6$.  Then
the $552$ points are the vectors $v\in \Lambda_{24}$ satisfying
$|v|^2 = |w-v|^2 = 4$.  These points all lie on a hyperplane, but
it does not pass through the origin, so subtract $w/2$ from each
vector to obtain points on a sphere centered at the origin.

The Higman-Sims configuration of $100$ points in $\R^{22}$ can be
constructed simi\-larly, following \cite{Wil}.  Choose $w_1$ and
$w_2$ in $\Lambda_{24}$ satisfying $|w_1|^2=|w_2|^2=6$ and
$|w_1-w_2|^2=4$. Then there are $100$ points $v \in \Lambda_{24}$
satisfying $|v|^2 = |w_1-v|^2 = |w_2-v|^2=4$, and they form the
Higman-Sims configuration (on an affine subspace as above).

One might imitate the last two constructions by using the $E_8$
lattice instead of the Leech lattice and replacing the norms $6$
and $4$ with the smallest two norms in $E_8$, namely $4$ and $2$.
That works, but it merely yields the cross polytope in $\R^7$ and
the simplex in $\R^6$.

The last line of the table describes a remarkable family of sharp
configurations from \cite{CGS}.  They are the only known sharp
configurations not derived from regular polytopes, the $E_8$
lattice, or the Leech lattice.  The parameter $q$ must be a prime
power. When $q=2$, this arrangement is the Schl\"afli configuration
(and it is a spherical $4$-design), but for $q>2$ it is different
from all the other entries in the table. Points in the
configuration correspond to totally isotropic $2$-dimensional
subspaces of a $4$-dimensional Hermitian space over $\F_{q^2}$,
with the distances between points determined by the dimensions of
the intersections of the corresponding subspaces (inner product
$-1/q$ corresponds to intersection dimension~$1$). This
construction generalizes to higher-dimensional Hermitian spaces,
but the $891$-point configuration in $\R^{22}$ appears to be the
only other sharp configuration that can be obtained in this way,
from totally isotropic $3$-dimensional subspaces of a
$6$-dimensional Hermitian space over $\F_4$.

Antipodal sharp configurations are the same as antipodal tight
spherical designs, by Theorem~6.8 in \cite{DGS}.  Much progress
has been made towards classifying such designs; see \cite{BMV} for
more details.  In particular, Table~\ref{tab:sharp} contains all
of the antipodal sharp configurations in at most $103$ dimensions.

Table~\ref{tab:sharp} coincides with the list of known cases in
which the linear programming bounds for spherical codes are sharp;
the list comes from \cite[p.~621]{L3}, except for the $600$-cell,
which was dealt with in \cite{A3} and \cite{E}.  We conjecture
that our techniques apply to a configuration if and only if the
linear programming bounds for spherical codes prove a sharp bound
for it.

In our main theorem, we prove that sharp arrangements and also the
$600$-cell are global potential energy minima for every completely
monotonic potential function:

\begin{theorem} \label{theorem:main}
Let $f \co (0,4] \to \R$ be completely monotonic, and let $\CodeS
\subset S^{n-1}$ be a sharp arrangement or the vertices of a
regular $600$-cell. If $\CodeT \subset S^{n-1}$ is any subset
satisfying $|\CodeT|=|\CodeS|$, then
\begin{equation}
\label{eq:bound} \sum_{x,y \in \CodeT, x \ne y} f\big(|x-y|^2\big)
\ge \sum_{x,y \in \CodeS, x \ne y} f\big(|x-y|^2\big).
\end{equation}
If $f$ is strictly completely monotonic, then equality in
\eqref{eq:bound} implies that $\CodeT$ is also a sharp
configuration or the vertices of a $600$-cell (whichever $\CodeS$
is) and the same distances occur in $\CodeS$ and $\CodeT$. In that
case, if $\CodeS$ is listed in Table~\ref{tab:sharp} but not on
the last line, then $\CodeT = A\CodeS$ for some $A \in O(n)$
(i.e., $\CodeT$ and $\CodeS$ are isometric).
\end{theorem}

Uniqueness does not necessarily hold if the hypothesis that $f$
must be strictly completely monotonic is removed (consider a
constant potential function).  Furthermore, the configurations
from the last line in Table~\ref{tab:sharp} are not always unique.
The construction in \cite{CGS} builds such a configuration given
any generalized quadrangle with parameters $(q,q^2)$.  For $q \le
3$ such quadrangles are unique (see Section~9 in \cite{CGS}), as
are the spherical codes derived from them, but for larger $q$ that
is not always true. Whenever distinct generalized quadrangles
exist, the construction yields distinct spherical codes.  In
particular, as explained in \cite{CGS}, one can construct such
quadrangles from ovoids in $\F_q\Proj^3$.  When $q=2^k$ with $k$
odd and $k>1$, there are at least two distinct ovoids, namely the
elliptic quadric ovoid and the Suzuki-Tits ovoid (see \cite{T}).
They yield distinct generalized quadrangles and thus distinct
universally optimal spherical codes.

To illustrate the techniques used to prove
Theorem~\ref{theorem:main}, in this section we will prove the
simplest possible case, namely when $|\CodeS| \le n+1$ (in which
case $\CodeS$ is a simplex). The full proof is in
Section~\ref{sec:proof}.

\begin{proof}[Proof of special case]
Let $N = |\CodeS|$, and suppose $N \le n+1$.  Then $\CodeS$ is a
regular simplex with inner product $-1/(N-1)$ between distinct
points, so the squared Euclidean distance between distinct points
is $2+2/(N-1)$ (because $|x-y|^2 = 2-2\langle x,y \rangle$ when
$|x|^2=|y|^2=1$).  For this special case, we only require $f$ to
be decreasing and convex, rather than completely monotonic.

Let
$$
h(x) = f\big(2+2/(N-1)\big) +
f'\big(2+2/(N-1)\big)\big(x-(2+2/(N-1))\big),
$$
so $h$ is the tangent line to $f$ at $2+2/(N-1)$.  Because $f$ is
convex, $f(x) \ge h(x)$ for all $x \in (0,4]$, and if $f$ is
strictly convex, then equality holds if and only if $x =
2+2/(N-1)$.

Suppose $\CodeT$ is any subset of $S^{n-1}$ with $|\CodeT|=N$.
Then
$$
\sum_{x,y \in \CodeT, x \ne y} f\big(|x-y|^2\big) \ge \sum_{x,y
\in \CodeT, x \ne y} h\big(|x-y|^2\big),
$$
and the right side equals
\begin{eqnarray*}
N(N-1)f\big(2+2/(N-1)\big)\\
\phantom{} -
(2+2/(N-1))N(N-1)f'\big(2+2/(N-1)\big)\\
\phantom{} + f'\big(2+2/(N-1)\big)\sum_{x,y \in \CodeT} |x-y|^2.
\end{eqnarray*}
We have
\begin{eqnarray*}
\sum_{x,y \in \CodeT} |x-y|^2 &=& \sum_{x,y \in \CodeT}
\big(2-2\langle x,y \rangle\big)\\
&=& 2N^2 - 2 \left|\sum_{x \in \CodeT} x\right|^2\\
&\le& 2N^2.
\end{eqnarray*}
Because $f'\big(2+2/(N-1)\big) \le 0$, we find that the potential
energy of $\CodeT$ is at least
\begin{eqnarray*}
N(N-1)f\big(2+2/(N-1)\big)\\
\phantom{} -
(2+2/(N-1))N(N-1)f'\big(2+2/(N-1)\big)\\
\phantom{} + 2N^2f'\big(2+2/(N-1)\big),
\end{eqnarray*}
which equals
$$
N(N-1)f\big(2+2/(N-1)\big).
$$

That is the potential energy of $\CodeS$, so we have proved that
$$
\sum_{x,y \in \CodeT, x \ne y} f\big(|x-y|^2\big) \ge \sum_{x,y
\in \CodeS, x \ne y} f\big(|x-y|^2\big).
$$
When $f$ is strictly convex, equality holds only if
$|x-y|^2=2+2/(N-1)$ for all $x,y \in \CodeT$ with $x \ne y$, in
which case $\CodeT$ is a regular simplex.
\end{proof}

We use the term \textit{universal optimality\/} to describe the
conclusion of Theorem~\ref{theorem:main}:

\begin{definition}
A finite subset $\CodeS \subset S^{n-1}$ is \textit{universally
optimal\/} if it (weakly) minimizes potential energy among all
configurations of $|\CodeS|$ points on $S^{n-1}$ for each
completely monotonic potential function.
\end{definition}

As the name suggests, every universally optimal spherical code is
indeed an optimal spherical code, in the sense of maximizing the
minimal distance between points. For the potential function $f(r)
= 1/r^s$, the leading term in the large-$s$ asymptotics of the
potential energy comes from the minimal distance. If there were a
spherical code with the same number of points but a larger minimal
distance, then it would have lower potential energy when $s$ is
sufficiently large.

As we mentioned above, the optimality of sharp configurations as
spherical codes is not new. Levenshtein proved in \cite{L2} that
every sharp configuration is a maximal spherical code, by using
linear programming bounds.  For the $E_8$ and Leech minimal
vectors, this amounts to the solution of the kissing problem by
Levenshtein and independently by Odlyzko and Sloane in
\cite{L,OS}; see \cite{L3} for an overview of the other cases,
except for the $600$-cell, for which optimality follows from the
bound in \cite{B} (see also \cite{A3} and \cite{E} for a solution
using linear programming bounds).  Our theorem can thus be viewed
as a continuum of optimality results connecting the kissing
problem and Levenshtein's theorem to the results of \cite{KY1} and
\cite{A2}.

Certain other known optimality results can also be derived.  For
example, \cite{KY1} and \cite{A2} prove that the $E_8$ and Leech
minimal vectors maximize the sum of the distances between pairs of
points. That is the case $f(r)=2-r^{1/2}$ in
Theorem~\ref{theorem:main}.  Similarly, taking $f(r) = \log(4/r)$
shows that sharp configurations maximize the product of the
distances between pairs of distinct points.

When testing whether a configuration is universally optimal, it is
helpful to be able to restrict to a smaller class of potential
functions.  Theorem~9b in \cite[p.~154]{W} implies that on each
compact subinterval of $(0,4]$, every completely monotonic
function on $(0,4]$ can be approximated uniformly by nonnegative
linear combinations of the functions $r \mapsto (4-r)^k$ with $k
\in \{0,1,2,\dots\}$.  For example,
$$
\frac{1}{r^s} = \sum_{k \ge 0}
\binom{s+k-1}{k}\frac{(4-r)^k}{4^{k+s}}.
$$
Thus, a configuration is universally optimal if and only if it is
optimal for each of the potential functions $f(r) = (4-r)^k$.

The conclusion of optimality under all completely monotonic
potential functions is quite strong, and computer experiments
suggest that it is rarely true in other cases.  For example, it
provably fails in the case of five points on $S^2$.  Sch\"utte and
van der Waerden proved in \cite{SW} that the minimal angle in such
a code cannot be greater than $\pi/2$. Their analysis shows that
all optimal codes consist of a pair of antipodal points with three
points on the orthogonal plane.  Among such codes, the only one
that could be universally optimal is the one in which the three
points form an equilateral triangle. However, when $f(r) =
(4-r)^7$, there is a configuration with lower energy: take a single
point together with four points equidistant from it and forming a
square in the opposite hemisphere (the size of the square is
chosen so as to minimize energy).  Thus, the only candidate for a
five-point universal optimum is not universally optimal.
Furthermore, it is not even a local minimum for energy when $f(r)
= (4-r)^k$ with $k$ sufficiently large or $f(r) = 1/r^s$ with $s$
sufficiently large.

Leech proved a far stronger theorem in \cite{Leech}.  He
classified the configurations on $S^2$ that are in equilibrium
under all force laws, which is a necessary but not sufficient
condition for universal optimality.  His theorem implies that
there are no universal optima in $\R^3$ beyond those in
Table~\ref{tab:sharp}, by ruling out all cases except for five
points (see \cite[p.~89]{Leech}). We conjecture that for each $n
\ge 3$ there are only finitely many universal optima in $\R^n$,
but we cannot prove it even for $n=4$.

We do not know of any other spherical codes that could be added to
Table~\ref{tab:sharp}, but some may well exist.  One surprising
case in which universal optimality fails is the vertices of the
regular $24$-cell on $S^3$ (equivalently, the $D_4$ root system).
It is not a sharp configuration, because it has four inner
products besides $1$ while it is only a spherical $5$-design, but
one might still hope that it would be universally optimal. It
appears to be optimal for many inverse power law potentials, but
even this weaker conclusion is surprisingly subtle: for example,
when one simulates $24$ randomly chosen points on $S^3$ under the
potential function $f(r)=1/r$, more than 90\% of the time one
reaches a local optimum with potential energy $334.096\dots$ (an
algebraic number of degree $20$),
which is only slightly worse than the potential energy $334$ of
the $24$-cell. In other words, there is a suboptimal arrangement
with a much larger basin of attraction.  This arrangement belongs
to an explicit one-parameter family of configurations analyzed in
\cite{CCEK}.  In that paper, some of them were shown to improve on
the $24$-cell for certain completely monotonic potential
functions, including $f(r) = (4-r)^k$ with $8 \le k \le 13$ and
$f(r) = e^{-3r}$. The $24$-cell is therefore not universally
optimal.

The remaining regular polytopes (the cubes in $\R^n$ for $n \ge
3$, the dodecahedron in $\R^3$, and the $120$-cell in $\R^4$) are
also not universally optimal, for a less interesting reason: they
are not even optimal spherical codes. Antiprisms improve on cubes,
and \cite{Sl} deals with the other two cases. Thus, a
full-dimensional regular polytope is universally optimal if and
only if all its faces are simplices.

One might also wonder about the $(4,10,1/6)$ code that is the
kissing configuration of the universally optimal $(5,16,1/5)$
code.  It consists of the midpoints of the edges of a regular
simplex. This code is not universally optimal, because for the
potential functions $f(r) = (4-r)^k$  with $3 \le k \le 6$,
arranging two regular pentagons on orthogonal planes yields lower
potential energy.  For $k \le 2$ their potential energies agree;
for $k \ge 7$, the $(4,10,1/6)$ code is superior to the two
pentagons, but we do not know whether it is optimal in these
cases.  Note that the kissing configuration of this $(4,10,1/6)$
code is a $(3,6,1/7)$ code that is not even optimal, let alone
universally optimal.

The $(4,10,1/6)$ code described above is a spherical $2$-design in
which only two inner products occur between distinct points.  It
follows that replacing the bound $2m-1$ with $2m-2$ in
Definition~\ref{def:sharp} would make Theorem~\ref{theorem:main}
false.

Probably the $(21,336,1/5)$ kissing configuration of the
$(22,891,1/4)$ code is not universally optimal, but in this case
we do not know a better configuration. We also do not know whether
universal optimality holds for the $(20,105,1/8)$ kissing
configuration of the $(21,162,1/7)$ code or the $(21,77,1/12)$
kissing configuration of the $(22,100,1/11)$ code.  These
configurations are not sharp, so Theorem~\ref{theorem:main} does
not apply.

We do not know whether there exist universally optimal spherical
codes whose universal optimality cannot be proved by our methods.
It seems likely that such codes exist, but we have no candidates
to propose.

Another important question not addressed by our methods is whether
there are other local minima for potential energy in the cases
listed in Table~\ref{tab:sharp}.  In certain special cases the
existence of other minima is trivial; for example, for $f(r)=4-r$,
a configuration minimizes potential energy if and only if its
centroid is the origin, so there is a positive-dimensional space
of global minima (when $N>3$).  In general, that is true for $f(r)
= (4-r)^k$ if $k$ is small relative to the number of points in the
configuration.  (It was proved in \cite{SZ} that spherical
$k$-designs exist on $S^{n-1}$ for all sufficiently large numbers
of points.  Spherical $k$-designs automatically minimize energy
for $f(r) = (4-r)^k$, and because the set of such designs is
closed under rotation and taking unions they exist in
positive-dimensional families once the number of points is large
enough.) For sufficiently large $k$ or for inverse power laws, the
set of local minima appears to be finite, but sometimes there are
several. When the number of points is large, there can be vast
numbers of local minima.  For example, experiments suggest that
for $120$ points on $S^3$ under the harmonic power law potential,
there are more than five thousand distinct local minima, and
perhaps far more than that.

It is natural to ask whether the requirement of complete
monotonicity could be weakened.  For example, regular simplices
are optimal for all decreasing, convex potential functions,
although not for all decreasing potential functions. That does not
always hold. Specifically, the icosahedron is not even a local
minimum for the decreasing, convex potential function $f(r) =
r^3-48r+128$, but rather a saddle point.  One can show using
Proposition~\ref{prop:lpbound} that the global minimum for twelve
points in $\R^3$ is obtained degenerately by letting them coincide
in triples at the vertices of a tetrahedron.  More generally,
repeating the vertices of a regular simplex in $\R^n$ is optimal
whenever the number of points is a multiple of $n+1$.  This
potential function therefore behaves quite differently from
completely monotonic potential functions.

The potential function $f(r) = r^3-48r+128$ from the previous
paragraph is convex, but $r \mapsto f(r^2)$ is not (so the
potential energy is not convex as a function of distance, as
opposed to distance squared). At the cost of slightly greater
complexity one can obtain that stronger property: if $f(r) = r^7 -
7\cdot2^{14}\cdot \sqrt{r} + 13\cdot2^{14}$, then $r \mapsto
f(r^2)$ is decreasing and convex, as is $f$, but the icosahedron
is not a local minimum for potential energy.

Although it is not enough to assume that the potential function is
decreasing and convex, in each case only finitely many derivative
conditions are needed.  For each sharp configuration, our methods
produce an explicit bound $K$ such that optimality holds for a
potential function $f \in C^K((0,4])$ whenever $(-1)^kf^{(k)}(r)
\ge 0$ for all $r \in (0,4]$ and $k \le K$. We will not discuss
this issue further, but it is worth recording the result: for a
sharp configuration with $m$ inner products between distinct
points, we can take $K = 2m$; furthermore, if $(-1)^k f^{(k)}(r)
> 0$ for all $r \in (0,4)$ and $k \le 2m$, then every
configuration of the same size and dimension with this potential
energy is sharp (and the same distances occur between points in it
as in the original configuration). For antipodal configurations a
more refined construction can be used to reduce both of these
bounds by $1$, but it is implausible that they could be
substantially reduced using our methods.

We conclude the introduction with one negative result in arbitrary
dimensions, namely that there are no universal optima in between
the simplex and cross polytope:

\begin{proposition}
If $n+1 < N < 2n$, then there is no $N$-point universally optimal
configuration on $S^{n-1}$.
\end{proposition}

\begin{proof}
An optimal code of this size has the following structure (see, for
example, Theorem~6.2.1 in \cite{BJr} and the remark following it).
There are $N-n$ pairwise orthogonal subspaces of $\R^n$ whose span
is $\R^n$, each containing $d+1$ points of the code if its
dimension is $d$. If the code is universally optimal, then each of
these subspaces must contain a regular simplex: because of
orthogonality, distances between points in different subspaces are
constant, so within each subspace we must have a universally
optimal code.  For the same reason, every union of some of the
component simplices must be universally optimal as well. Thus,
without loss of generality we can assume that $N=n+2$. We can also
assume that $n \ge 4$ because the case of $n=3$ and $N=5$ was
dealt with earlier in this section.

Suppose the two regular simplices have $i+1$ and $n-i+1$ points.
For the potential function $f(r) = (4-r)^k$, the energy of the
configuration equals
$$
i(i+1)\left(2-\frac{2}{i}\right)^k+
(n-i)(n-i+1)\left(2-\frac{2}{n-i}\right)^k+(i+1)(n-i+1)2^{k+1}.
$$
For $k=2$ that equals $4n(n+1)-4n/\big(i(n-i)\big)$, from which it
follows that the potential energy is minimized exactly when $i =
\lfloor n/2 \rfloor$ or $i = \lceil n/2 \rceil$.  As $k \to
\infty$, the energy is asymptotic to $(i+1)(n-i+1)2^{k+1}$, which
is minimized exactly when $i=1$ or $i=n-1$ (and hence not at $i =
\lfloor n/2 \rfloor$ or $i = \lceil n/2 \rceil$, because $n \ge
4$). Thus, universal optimality fails, as asserted.
\end{proof}

\section{Background}

In this section we briefly review the tools required for the proof
of Theorem~\ref{theorem:main}.

\subsection{Hermite interpolation}

Recall that \textit{Hermite interpolation\/} is a generalization
of Lagrange interpolation in which one computes a polynomial that
matches not only the values but also some derivatives of a
function at some specified points.  More precisely, suppose we are
given $f \in C^\infty([a,b])$, distinct points $t_1,\dots,t_m \in
[a,b]$, and positive integers $k_1,\dots,k_m$.  We wish to find a
polynomial of degree less than $D=k_1+k_2+\dots+k_m$ such that for
$1 \le i \le m$ and $0 \le k < k_i$,
$$
p^{(k)}(t_i) = f^{(k)}(t_i).
$$
In other words, $p$ and $f$ agree to order $k_i$ at $t_i$. A
unique such polynomial always exists, because the linear map that
takes a polynomial $p$ of degree less than $D$ to
$$
\big(p(t_1),p'(t_1),\dots,p^{(k_1-1)}(t_1),p(t_2),p'(t_2),
\dots,p^{(k_m-1)}(t_m)\big)
$$
is injective (and hence surjective as well).

One fundamental fact about Hermite interpolation is the following
remainder formula (see Theorem~3.5.1 in \cite[p.~68]{D}):

\begin{lemma}\label{lemma:remainder}
Under the hypotheses above, for each $t \in [a,b]$ there exists
$\xi \in (a,b)$ such that $\min(t,t_1,\dots,t_m) < \xi <
\max(t,t_1,\dots,t_m)$ and
$$
f(t) - p(t) = \frac{f^{(D)}(\xi)}{D!}(t-t_1)^{k_1}\dots
(t-t_m)^{k_m}.
$$
\end{lemma}

Because the standard proof is short, we reproduce it here for the
convenience of the reader:

\begin{proof}
Let
$$
g(t) = \frac{f(t) - p(t)}{(t-t_1)^{k_1}\dots (t-t_m)^{k_m}}.
$$
For fixed $t \in [a,b] \setminus \{t_1,\dots,t_m\}$ (the lemma is
trivial for $t \in \{t_1,\dots,t_m\}$), consider the function
$$
s \mapsto f(s) - p(s) - g(t)(s-t_1)^{k_1}\dots (s-t_m)^{k_m}.
$$
By construction, for each $i$ it vanishes at $t_i$ to order $k_i$,
and it also vanishes at $t$, so it has $D+1$ roots in the interval
$$
[\min(t,t_1,\dots,t_m),\max(t,t_1,\dots,t_m)].
$$
By iterated use of Rolle's theorem, there exists $\xi$ in the
interior of this interval at which the $D$-th derivative vanishes.
At that point,
$$
f^{(D)}(\xi) - g(t) D! = 0,
$$
which gives the desired formula for $g(t)$.
\end{proof}

This formula will play a fundamental role later in the paper, but
we will also need a stronger result.  Recall that a function is
\textit{absolutely monotonic\/} on an interval if it is $C^\infty$
and it and all its derivatives are nonnegative on that interval
(see \cite[p.~144]{W}). Unfortunately, the terms ``completely
monotonic'' and ``absolutely monotonic'' are easy to confuse, but
they are standard terms.  As a mnemonic, associate ``absolutely
monotonic'' with ``absolute value'' to remember which one has no
signs.  We call a function \textit{strictly absolutely
monotonic\/} on an interval if it and all its derivatives are
nonnegative on the interval and strictly positive on its interior.

\begin{proposition}\label{prop:absmonred}
Under the hypotheses above, suppose further that $f$ is absolutely
monotonic on $(a,b)$.  Then
$$
\frac{f(t)-p(t)}{(t-t_1)^{k_1}\dots (t-t_m)^{k_m}}
$$
is also absolutely monotonic on $(a,b)$.  If $f$ is strictly
absolutely monotonic, then
$$
\frac{f(t)-p(t)}{(t-t_1)^{k_1}\dots (t-t_m)^{k_m}}
$$
is as well.
\end{proposition}

Of course, when we write quotients such as this one, we define
them by continuity when $t = t_i$ for some $i$ so that they become
$C^\infty$ functions. It seems plausible that
Proposition~\ref{prop:absmonred} is known, but we do not know
where to find it in the literature. The hypothesis that
$t_1,\dots,t_m$ lie in the interval $[a,b]$ is crucial (it is part
of the setup above, but it is worth reiterating it here).

Before proving Proposition~\ref{prop:absmonred}, it is convenient
to introduce systematic notation for Hermite interpolations. Given
a polynomial $g$ with $\deg(g) \ge 1$, let $H(f,g)$ denote the
polynomial of degree less than $\deg(g)$ that agrees with $f$ at
each root of $g$ to the order of that root.  (If $f$ is a
polynomial, then $H(f,g)$ is its remainder modulo $g$ in
polynomial division.) In other words, if $g(t)$ vanishes to order
$k$ at $t=s$, then
$$
f(t) - H(f,g)(t) = O\big((t-s)^{k}\big)
$$
as $t \to s$.  It follows that the function $Q(f,g)$ defined by
$$
Q(f,g)(t) = \frac{f(t)-H(f,g)(t)}{g(t)}
$$
extends to a $C^\infty$ function at the roots of $g$.
Proposition~\ref{prop:absmonred} is the assertion that if
$g(t)=(t-t_1)^{k_1}\dots (t-t_m)^{k_m}$, then $Q(f,g)$ is
absolutely monotonic.

\begin{proof}
The proof rests on two key properties. The first is that
$$
Q(f,g_1g_2) = Q(Q(f,g_1),g_2),
$$
which follows from the uniqueness of the Hermite interpolation:
$$
Q(Q(f,g_1),g_2) = \frac{f-H(f,g_1)-g_1H(Q(f,g_1),g_2)}{g_1g_2},
$$
but $H(f,g_1g_2)$ is the unique polynomial of degree less than
$\deg(g_1)+\deg(g_2)$ such that $(f-H(f,g_1g_2))/(g_1g_2)$ extends
to a $C^\infty$ function everywhere.  Thus, $H(f,g_1g_2)
=H(f,g_1)+g_1H(Q(f,g_1),g_2)$ and $Q(f,g_1g_2) = Q(Q(f,g_1),g_2)$.

The second property is that if $g_0(t) = (t-s_0)^n$, then
$$
Q(f,g_0)(s_0) = \frac{f^{(n)}(s_0)}{n!}.
$$
This follows immediately from the Taylor expansion.

We must show that for $n \ge 0$ and $t \in (a,b)$,
$$
Q(f,g)^{(n)}(t) \ge 0.
$$
When $n=0$, that follows from Lemma~\ref{lemma:remainder}, which
asserts the existence of $\xi \in (a,b)$ (depending on $t$) such
that
$$
Q(f,g)(t) = \frac{f^{(\deg(g))}(\xi)}{\deg(g)!} \ge 0.
$$
For larger $n$, given $s_0 \in (a,b)$, define $g_0(t) = (t-s_0)^n$.
Then
$$
\frac{Q(f,g)^{(n)}(s_0)}{n!} = Q(Q(f,g),g_0)(s_0) = Q(f,gg_0)(s_0)
\ge 0.
$$
This completes the proof when $f$ is absolutely monotonic. The
case of strictly absolutely monotonic functions works exactly the
same way.
\end{proof}

It follows immediately from continuity that the statement of
Proposition~\ref{prop:absmonred} also holds when the open interval
$(a,b)$ is replaced with a half-open or closed interval.

\subsection{Positive-definite kernels}
\label{subsec:pdk}

A continuous function $K \co S^{n-1} \times S^{n-1} \to \R$ is a
\textit{positive-definite kernel\/} if for all $k$ and all
$x_1,\dots,x_k \in S^{n-1}$, the $k \times k$ matrix whose $i,j$
entry is $K(x_i,x_j )$ is positive semidefinite.  In other words,
for all $t_1,\dots,t_k \in \R$,
$$
\sum_{1 \le i,j \le k} t_i t_j K(x_i,x_j) \ge 0.
$$
The most important special case is when $t_1 = \dots = t_k = 1$:
for every finite subset $\CodeS \subset S^{n-1}$,
$$
\sum_{x,y \in \CodeS} K(x,y) \ge 0.
$$

Positive-definite kernels $K$ such that $K(x,y)$ depends only on
the distance $|x-y|$ between $x$ and $y$ will be particularly
important. There is a simple representation-theoretic construction
of such kernels, which Schoenberg proved gives all of them (see
\cite{S}).  As a unitary representation of $O(n)$, the Hilbert
space $L^2(S^{n-1})$ breaks up as a completed orthogonal direct
sum of infinitely many finite-dimensional representations:
$$
L^2(S^{n-1}) = \widehat{\bigoplus_{\ell \ge 0}} \,\, V_\ell,
$$
where $V_\ell$ is the space of spherical harmonics of degree
$\ell$ (i.e., restrictions to $S^{n-1}$ of homogeneous polynomials
on $\R^n$ of degree $\ell$ that are in the kernel of the Euclidean
Laplacian). Note that spherical harmonics of different degrees are
orthogonal, which will play a useful role later.  It is also
important to know that every polynomial on $\R^n$ has the same
restriction to $S^{n-1}$ as some unique linear combination of
spherical harmonics. For more details on spherical harmonics, see
\S{}2 of Chapter~IV in \cite{StW}.

For each $\ell \ge 0$ and each $x \in S^{n-1}$, there is a unique
\textit{reproducing kernel\/} $\ev_{\ell,x} \in V_\ell$, defined
by requiring that for all $f \in V_\ell$,
$$
\langle f, \ev_{\ell,x}\rangle = f(x).
$$
(We use Hermitian forms that are conjugate-linear in the second
variable.)  To construct $\ev_{\ell,x}$ explicitly, choose an
orthonormal basis $e_{\ell,1},\dots,e_{\ell,\dim V_\ell}$ of
$V_\ell$ and define
$$
\ev_{\ell,x}(y) = \overline{e_{\ell,1}(x)}e_{\ell,1}(y) + \dots +
\overline{e_{\ell,\dim V_\ell}(x)}e_{\ell,\dim V_\ell}(y).
$$

Define the function $K_\ell \co S^{n-1} \times S^{n-1} \to \C$ by
$$
K_\ell(x,y) = \langle \ev_{\ell,x}, \ev_{\ell,y} \rangle =
\ev_{\ell,x}(y).
$$
It is straightforward to check that $\ev_{\ell,Ax} = \ev_{\ell,x}
\circ A^{-1}$ for $A \in O(n)$.  Because $O(n)$ acts
distance-transitively on $S^{n-1}$ (i.e., if $|u-v|=|x-y|$, then
there exists $A \in O(n)$ such that $Au=x$ and $Av=y$),
$K_\ell(x,y)$ depends only on the distance between $x$ and $y$.
That implies that it is real-valued, since
$$
K_\ell(x,y) = \overline{K_\ell(y,x)}
$$
follows immediately from the definition while $K_\ell(x,y) =
K_\ell(y,x)$ since $|x-y|=|y-x|$.

The key property of $K_\ell$ is that it is a positive-definite
kernel:
$$
\sum_{1 \le i,j \le k} t_i t_j K_\ell(x_i,x_j) = \left|\sum_{1 \le
i \le k} t_i \ev_{\ell,x_i}\right|^2 \ge 0.
$$
All convergent nonnegative infinite linear combinations of these
functions are still positive-definite kernels, and Schoenberg
proved that those are the only continuous positive-definite
kernels that depend only on the distance between points.

One can compute $K_\ell$ explicitly.  Because $K_\ell(x,y)$
depends only on the distance between $x$ and $y$, we can write
$K_\ell(x,y) = C_\ell(\langle x,y \rangle)$ for some function
$C_\ell$.  To compute $C_\ell$, we use the orthogonality between
spherical harmonics of different degrees. If we fix a point on the
sphere and project orthogonally onto the line through the origin
and that point, then the surface measure on the sphere projects to
a constant times the measure
$$
(1-t^2)^{(n-3)/2} \, dt
$$
on the interval $[-1,1]$, and the orthogonality property amounts
to saying that $C_0,C_1,\dots$ are orthogonal polynomials with
respect to that measure, with $C_i$ having degree $i$. That
determines these functions up to multiplication by a positive
scalar (positive because $C_i(1) > 0$), and pinning down that
factor is irrelevant for our purposes.

These orthogonal polynomials are known as the ultraspherical (or
Gegenbauer) polynomials.  The standard notation is
$C_i^\lambda(t)$ for the degree $i$ polynomial orthogonal with
respect to $(1-t^2)^{\lambda-1/2}\, dt$ on $[-1,1]$ (for
$S^{n-1}$, $\lambda = n/2-1$). These polynomials are normalized so
that $C_0^\lambda(t)=1$ and $C_1^\lambda(t) = 2\lambda t$, and
they satisfy the recurrence
$$
i C^\lambda_i(t) = 2(i+\lambda-1)t C^\lambda_{i-1}(t) -
(i+2\lambda-2)C^\lambda_{i-2}(t)
$$
for $i \ge 2$.  See \cite[p.~303]{AAR} for more details on
ultraspherical polynomials, as well as Chapter~9 of \cite{AAR} for
spherical harmonics.  We always assume $\lambda \ge 0$
(equivalently, $n \ge 2$), because the one-dimensional case is
trivial and its ultraspherical polynomials are poorly normalized:
above degree $0$ they have negative leading coefficients, and
$C_1^{-1/2}$ is not a positive-definite kernel.

Whenever we refer to the \textit{ultraspherical coefficients\/} of
a function, we mean its coefficients in terms of ultraspherical
polynomials (where the parameter $\lambda$ is implicit in the
context).  A \textit{positive-definite function\/} will be a
function whose ultraspherical coefficients are nonnegative.

One fundamental fact about positive-definite functions is that
they are closed under taking products.  Equivalently, the product
of two ultraspherical polynomials is a nonnegative linear
combination of ultraspherical polynomials (for $\lambda \ge 0$).
By orthogonality, proving that amounts to showing that
\begin{equation} \label{eq:prodclosed}
\int_{-1}^1 C_i(t) C_j(t) C_k(t) \, (1-t^2)^{(n-3)/2} dt \ge 0
\end{equation}
for all $i,j,k \ge 0$. Let $\mu$ denote the surface measure on
$S^{n-1}$. Then using the expansion
$$
K_\ell(x,y) = \sum_{m=1}^{\dim V_\ell} \overline{e_{\ell,m}(x)}
e_{\ell,m}(y)
$$
in terms of a choice of orthonormal basis
$e_{\ell,1},\dots,e_{\ell,\dim V_\ell}$ for each $V_\ell$ shows
that \eqref{eq:prodclosed} holds if and only if
$$
\sum_{a=1}^{\dim V_i} \sum_{b=1}^{\dim V_j} \sum_{c=1}^{\dim V_k}
\int e_{i,a}(x) e_{j,b}(x) e_{k,c}(x) \, d\mu(x) \
\overline{e_{i,a}(y)} \, \overline{e_{j,b}(y)} \,
\overline{e_{k,c}(y)} \ge 0
$$
for some (equivalently, all) $y \in S^{n-1}$.  Integrating over
$y$ yields
$$
\sum_{a=1}^{\dim V_i} \sum_{b=1}^{\dim V_j} \sum_{c=1}^{\dim V_k}
\left|\int e_{i,a}(x) e_{j,b}(x) e_{k,c}(x) \, d\mu(x)
\right|^2\!\!,
$$
which is visibly nonnegative.

This result also follows from the Schur product theorem: the set
of positive semidefinite matrices of a given size is closed under
taking Hadamard products (i.e., taking the product entry by
entry). See Theorem~7.5.3 in \cite[p.~458]{HJ}.  That approach
leads to the proof of Schoenberg's characterization of
positive-definite kernels.  Suppose $(x,y) \mapsto f(\langle
x,y\rangle)$ is a positive-definite kernel.  By the Schur product
theorem, the same is true for $(x,y) \mapsto f(\langle x,y\rangle)
C_i(\langle x,y\rangle)$ for each $i$.  It follows that the
integral of this function over $(x,y) \in S^{n-1} \times S^{n-1}$
is nonnegative, and hence for fixed $y$ the integral over $x \in
S^{n-1}$ is nonnegative (because it is independent of $y$).  In
other words, the $i$-th ultraspherical coefficient of $f$ must be
nonnegative.

We call a polynomial \textit{strictly positive definite\/} if all
its ultraspherical coefficients are strictly positive (up to its
degree, of course).  The product of strictly positive-definite
polynomials is strictly positive definite;  the proof rests on the
fact that in the expansion of $C_i^\lambda C_j^\lambda$, the
coefficient of $C_{i+j}^\lambda$ is positive (because all these
polynomials have positive leading coefficients).

One useful application of ultraspherical polynomials is the
following simple test for spherical design strength (Theorem~5.5
in \cite{DGS}): a configuration $\CodeS \subset S^{n-1}$ is a
spherical $M$-design if and only if
\begin{equation} \label{eq:zerosum}
\sum_{x,y \in \CodeS} C_i^{n/2-1}(\langle x,y \rangle) = 0
\end{equation}
for $1 \le i \le M$.  That property depends only on the distance
distribution of $\CodeS$ (i.e, which distances occur between
points in $\CodeS$ and the multiplicities of those distances). To
prove this characterization of designs, observe that the identity
$$
\sum_{x,y \in \CodeS} C_i^{n/2-1}(\langle x,y \rangle) = \sum_{x,y
\in \CodeS} \langle \ev_{i,x},\ev_{i,y}\rangle = \left|\sum_{x \in
\CodeS} \ev_{i,x}\right|^2
$$
shows that \eqref{eq:zerosum} holds if and only if all of the
degree $i$ spherical harmonics sum to zero over $\CodeS$.

\section{Orthogonal polynomials}

The proof of Theorem~\ref{theorem:main} involves expansions as
nonnegative linear combinations of orthogonal polynomials.  In
this section we prove a theorem about such expansions that will
play a fundamental role later in the paper.  We prove it in
greater generality than is required for
Theorem~\ref{theorem:main}, because it is no more difficult and
the proof is of independent interest.

Let $\mu$ be any Borel measure on $\R$ such that all polynomials
are integrable with respect to $\mu$ and for all polynomials $p$,
$$
\int p(t)^2 \, d\mu(t) > 0
$$
if $p$ is not identically zero.  Equivalently, for all polynomials
$p$ such that $p(t) \ge 0$ for all $t$ but $p$ is not identically
zero,
$$
\int p(t) \, d\mu(t) > 0.
$$
(Every such polynomial is a sum of two squares of polynomials; see
for example Problem~44 in Part~VI of \cite{PS} and its solution.)

Let $p_0,p_1,\dots$ be the monic orthogonal polynomials for $\mu$,
with $\deg(p_i)=i$.  In other words, for $i \ne j$
$$
\int p_i(t) p_j(t) \, d\mu(t) = 0.
$$
It is known that for each $\alpha \in \R$, the polynomial $p_n +
\alpha p_{n-1}$ has $n$ distinct real roots, which are interlaced
with the roots of $p_{n-1}$ (see Theorem~3.3.4 in \cite{Sz}).

\begin{theorem} \label{theorem:orthog}
Let $\alpha$ be any real number, and let
$$
r_1 < r_2 < \dots < r_n
$$
be the roots of $p_n + \alpha p_{n-1}$.  Then for $k < n$,
$$
\prod_{i=1}^k (t-r_i)
$$
has positive coefficients in terms of
$p_0(t),p_1(t),\dots,p_k(t)$.
\end{theorem}

The proof is conceptually simple, but it involves some
technicalities.  We begin with the simplest case, namely $k=n-1$:

\begin{proposition} \label{prop:basecase}
Let $\alpha$ be any real number, and let $r$ be the largest root
of $p_n + \alpha p_{n-1}$.  Then
$$
\frac{p_n(t) + \alpha p_{n-1}(t)}{t-r}
$$
has positive coefficients in terms of
$p_0(t),p_1(t),\dots,p_{n-1}(t)$.
\end{proposition}

Proposition~\ref{prop:basecase} follows from the
Christoffel-Darboux formula (Theorem~5.2.4 in \cite{AAR}), but we
prove it directly to illustrate the technique we will apply to the
general case.

\begin{proof}
Define $c_0,\dots,c_{n-1}$ so that
$$
\frac{p_n(t) + \alpha p_{n-1}(t)}{t-r} = \sum_{\ell=0}^{n-1}
c_\ell p_\ell(t).
$$

For $\ell \le n-1$, it follows from orthogonality that
$$
\int (p_n(t) + \alpha p_{n-1}(t)) \frac{p_\ell(t) - p_\ell(r)}{t -
r} \, d\mu(t) = 0,
$$
because
$$
\frac{p_\ell(t) - p_\ell(r)}{t - r}
$$
is a polynomial of degree $\ell-1$ and hence orthogonal to both
$p_n(t)$ and $p_{n-1}(t)$. Thus,
$$
\int \frac{p_n(t) + \alpha p_{n-1}(t)}{t - r}p_\ell(t) \, d\mu(t)
= p_\ell(r) \int \frac{p_n(t) + \alpha p_{n-1}(t)}{t - r} \,
d\mu(t).
$$
In terms of the coefficients $c_0,\dots,c_{n-1}$,
$$
c_\ell \int p_\ell(t)^2 \, d\mu(t) = c_0 p_\ell(r)  \int d\mu(t).
$$

Because the roots of $p_n+\alpha p_{n-1}$ and $p_{n-1}$ are
interlaced, $r$ is greater than the largest root of $p_{n-1}$ and
hence greater than the largest root of $p_\ell$. It follows that
$p_\ell(r)>0$, and thus $c_0,\dots,c_{n-1}$ all have the same
sign. Comparing leading coefficients shows that $c_{n-1}=1$, so
every coefficient must be positive.
\end{proof}

There is a natural generalization of this proof to the case of
arbitrary $k$ in Theorem~\ref{theorem:orthog}.  The key idea is to
use orthogonality with respect to a signed measure, because for
$k<n$,
$$
\prod_{i=1}^k (t-r_i)
$$
is the monic orthogonal polynomial of degree $k$ for the signed
measure
$$
(t-r_{k+1})\dots(t-r_n) \, d\mu(t),
$$
assuming there is a unique such polynomial, as we prove below in
Lemma~\ref{lemma:signed} and Lemma~\ref{lemma:posdefbd}. To verify
that, we simply need to show that it is orthogonal to all
polynomials of degree less than $k$, which is equivalent with the
orthogonality of $p_n(t)+\alpha p_{n-1}(t)$ to all such
polynomials with respect to $d\mu(t)$. (For $k=n$ the conclusion
is false unless $\alpha=0$.)  More generally, we have the
following lemma:

\begin{lemma} \label{lemma:signedmeasure}
Suppose $\nu$ is a signed measure that has monic orthogonal
polynomials $q_0,q_1,\dots,q_{N+1}$, where $\deg(q_i) = i$. For $i
\le N$, if $q_i(r) \ne 0$ and $(t-r) \, d\nu(t)$ has a unique
degree $i$ monic orthogonal polynomial, then it equals
$$
\frac{q_{i+1}(t) + \alpha_i q_i(t)}{t-r},
$$
where $\alpha_i$ is chosen so that $q_{i+1}(r) + \alpha_i q_i(r) =
0$.
\end{lemma}

\begin{proof}
The orthogonality of this polynomial with respect to $(t-r)\,
d\nu(t)$ amounts to the orthogonality with respect to $d\nu(t)$ of
$q_{i+1}(t) + \alpha_i q_i(t)$ to all polynomials of degree at
most $i-1$.
\end{proof}

One might hope to imitate the proof of
Proposition~\ref{prop:basecase} as one removes successive roots of
$p_n+\alpha p_{n-1}$. In fact, this approach works, but it is not
obvious that it works, because of a fundamental technical
difficulty. As soon as one introduces linear factors such as $t-r$
into the measure, it is no longer a positive measure but rather a
signed measure.  Many of the basic properties of orthogonal
polynomials fail in the case of a signed measure.  Monic
orthogonal polynomials may not exist or be unique, and the roots
of $p_i$ and $p_{i+1}$ may not be interlaced. In addition, there
is an even more dramatic problem: the proof of
Proposition~\ref{prop:basecase} depends on the positivity of
$$
\int p_\ell(t)^2 \, d\mu(t)
$$
and
$$
\int d\mu(t),
$$
neither of which is generally positive for a signed measure.

Fortunately, the measures that actually arise in the proof will be
considerably better behaved than a typical signed measure.  In
particular, they will all have the following property for some
$N$:

\begin{definition}
Let $\nu$ be a signed Borel measure on $\R$ with respect to which
all polynomials are integrable. We say $\nu$ is \textit{positive
definite up to degree $N$} if for all polynomials $f$ such that
$\deg(f) \le N$,
$$
\int f(t)^2 \, d\nu(t) \ge 0,
$$
with equality only if $f$ is identically zero.
\end{definition}

If $\nu$ is positive definite up to degree $N$, then many of the
usual properties hold for orthogonal polynomials of degrees up to
$N+1$. Lemma~\ref{lemma:signed} covers the properties we will
require.

\begin{lemma}\label{lemma:signed}
Let $\nu$ be positive definite up to degree $N$.  Then there are
unique monic polynomials $q_0,q_1,\dots,q_{N+1}$ such that
$\deg(q_i)=i$ for each $i$ and
$$
\int q_i(t) q_j(t) \, d\nu(t) = 0
$$
for $i\ne j$.  For each $i$, $q_i$ has $i$ distinct real roots,
and the roots of $q_i$ and $q_{i-1}$ are interlaced.
\end{lemma}

The standard proofs for orthogonal polynomials from Chapter~5 of
\cite{AAR} can be adapted, if one uses the trick of writing a
nonnegative polynomial as a sum of two squares.  Here we use
instead a particularly elegant proof from \cite[p.~14]{Si}:

\begin{proof}
The orthogonal polynomials for $\nu$ are obtained by applying
Gram-Schmidt orthogonalization to $1,t,\dots,t^{N+1}$:
$$
q_i(t) = t^i - \sum_{j=0}^{i-1} q_j(t) \frac{\int s^i q_j(s) \,
d\nu(s)}{\int q_j(s)^2 \, d\nu(s)},
$$
which makes sense for $i \le N+1$ because the norm of $q_j$ in the
denominator is always nonzero.  They are unique because the
quotients of integrals are the only coefficients in such an
expansion that yield orthogonality.

The next step is to prove the three-term recurrence relation.  The
polynomial $q_i(t) - tq_{i-1}(t)$ is orthogonal to all polynomials
of degree less than $i-2$, and hence there exist constants $a_i$
and $b_i$ such that
$$
q_i(t) = (t+a_i)q_{i-1}(t) + b_i q_{i-2}(t).
$$
Multiplying by $q_{i-2}(t)$ and integrating yields
$$
b_i \int q_{i-2}(t)^2 \, d\nu(t) = -\int q_{i-1}(t) t q_{i-2}(t)
\, d\nu(t),
$$
by orthogonality.  The polynomials $tq_{i-2}(t)$ and $q_{i-1}(t)$
differ by a polynomial of degree at most $i-2$, and hence
$$
\int q_{i-1}(t) t q_{i-2}(t) \, d\nu(t) = \int q_{i-1}(t)^2 \,
d\nu(t).
$$
Combining the last two equations shows that $b_i < 0$.

We now prove by induction that the roots of $q_i$ and $q_{i-1}$
are real and interlaced.  The base case of $i=1$ is trivial.  For
the induction step, suppose the roots of $q_{i-1}$ and $q_{i-2}$
are real and interlaced.  For each root $r$ of $q_{i-1}$,
$$
q_i(r) = b_i q_{i-2}(r),
$$
by the three-term recurrence relation.  By hypothesis, $q_{i-2}$
alternates in sign at the roots of $q_{i-1}$.  Because $b_i$ is
negative, $q_i$ and $q_{i-2}$ have opposite signs at the roots of
$q_{i-1}$. That implies that $q_i$ must have a root between each
pair of consecutive roots of $q_{i-1}$.  Because $q_i(t)$ and
$q_{i-2}(t)$ have the same sign when $|t|$ is sufficiently large,
$q_i$ must also have a root greater than every root of $q_{i-1}$
and a root less than every root of $q_{i-1}$.  All $i$ roots of
$q_i$ are now accounted for, so all the roots of $q_i$ are real
and interlaced with those of $q_{i-1}$.  That completes the proof.
\end{proof}

The following minor variant of Gauss-Jacobi quadrature (see
Theorem 5.3.2 in \cite{AAR}) will prove crucial for understanding
the signed measures that arise in our proof.

\begin{lemma} \label{lemma:quadrature}
Let $\alpha$ be any real number, and let
$$
r_1 < r_2 < \dots < r_n
$$
be the roots of $p_n + \alpha p_{n-1}$.  Then there are positive
numbers $\lambda_1,\dots,\lambda_n$ such that for every polynomial
$f$ of degree at most $2n-2$,
$$
\int f(t) \,d\mu(t) = \sum_{i=1}^n \lambda_i f(r_i).
$$
\end{lemma}

\begin{proof}
There exist coefficients $\lambda_1,\dots,\lambda_n$ such that
$$
\int f(t) \,d\mu(t) = \sum_{i=1}^n \lambda_i f(r_i)
$$
whenever $\deg(f)<n$, because a polynomial of degree less than $n$
is completely determined by its values at $r_1,\dots,r_n$ and the
map taking these values to the integral is linear.

For $\deg(f) \le 2n-2$, there exist polynomials $g$ and $h$
satisfying 
$$f = (p_n+\alpha p_{n-1})g + h$$ with $\deg(g) \le n-2$
and $\deg(h) < n$.  By orthogonality,
$$
\int f(t) \,d\mu(t) = \int h(t) \,d\mu(t),
$$
and $f(r_i) = h(r_i)$ for each $i$ because $p_n(r_i) + \alpha
p_{n-1}(r_i)=0$ by definition.  It follows that
$$
\int f(t) \,d\mu(t) = \sum_{i=1}^n \lambda_i f(r_i)
$$
holds whenever $\deg(f) \le 2n-2$.

All that remains is to prove positivity for the coefficients.  For
fixed $i$, let
$$
f(t) = \prod_{j \ne i} (t-r_j)^2.
$$
Then
$$
\int f(t) \,d\mu(t) = \lambda_i f(r_i)^2,
$$
from which it follows that $\lambda_i>0$.
\end{proof}

Let $r_1 < \dots < r_n$ be the roots of $p_n + \alpha p_{n-1}$,
and for $0 \le j \le n$ define the measure $\mu_j$ by
$$
d\mu_j(t) = \prod_{i=0}^{j-1} (r_{n-i}-t) \, d\mu(t)
$$
(of course $\mu_0=\mu$).

\begin{lemma} \label{lemma:posdefbd}
For $0 \le j \le n-1$, the measure $\mu_j$ is positive definite up
to degree $n-j-1$.
\end{lemma}

\begin{proof}
If $\deg(f) \le n-1-j/2$, then by Lemma~\ref{lemma:quadrature},
$$
\int f(t)^2 \prod_{i=0}^{j-1} (r_{n-i}-t) \, d\mu(t) =
\sum_{i=1}^n \lambda_i f(r_i)^2 \prod_{\ell=0}^{j-1}
(r_{n-\ell}-r_i) \ge 0;
$$
equality holds if and only if $f$ vanishes at $r_1,\dots,r_{n-j}$,
which is impossible for a polynomial of degree at most $n-j-1$
unless it vanishes identically.
\end{proof}

Let $q_{j,i}$ be the monic orthogonal polynomial of degree $i$ for
$\mu_j$.  (Note that $q_{0,i} = p_i$.)  The polynomial $q_{j,i}$
is simplest when $i=n-j$, in which case we have
$$
q_{j,n-j}(t) = (t-r_1) \dots (t-r_{n-j}),
$$
as pointed out earlier. That implies that for $i < n-j$, the
largest root of $q_{j,i}$ is less than $r_{n-j}$, because the
roots are interlaced.

Equivalently, for $i \le n-j$, the largest root of $q_{j-1,i}$ is
less than $r_{n-j+1}$, so $q_{j-1,i}(r_{n-j+1}) \ne 0$.  Hence
Lemma~\ref{lemma:signedmeasure} implies that for $i \le n-j$,
there are constants $\alpha_{j,i}$ such that
$$
q_{j,i}(t) = \frac{q_{j-1,i+1}(t) + \alpha_{j,i}
q_{j-1,i}(t)}{t-r_{n-j+1}}.
$$

\begin{lemma} \label{lemma:induction}
For $1 \le j \le n$ and $i \le n-j$, the polynomial $q_{j,i}$ is a
positive linear combination of $q_{j-1,0},\dots,q_{j-1,i}$.
\end{lemma}

\begin{proof}
Define $c_0,\dots,c_i$ so that
$$
q_{j,i}(t) = \sum_{\ell = 0}^i c_\ell q_{j-1,\ell}(t).
$$
We argue as in the proof of Proposition~\ref{prop:basecase}. For
$\ell \le i$,
\begin{equation} \label{eq:ortho}
\int \big(q_{j-1,i+1}(t) + \alpha_{j,i} q_{j-1,i}(t)\big)
\frac{q_{j-1,\ell}(t)-q_{j-1,\ell}(r_{n-j+1})}{t-r_{n-j+1}} \,
d\mu_{j-1}(t) = 0
\end{equation}
by orthogonality, because
$$
\frac{q_{j-1,\ell}(t)-q_{j-1,\ell}(r_{n-j+1})}{t-r_{n-j+1}}
$$
is a polynomial of degree $\ell-1$, which is less than $i$.  It
now follows from
$$
q_{j,i}(t) = \frac{q_{j-1,i+1}(t) + \alpha_{j,i}
q_{j-1,i}(t)}{t-r_{n-j+1}}
$$
and \eqref{eq:ortho} that
$$
\int q_{j,i}(t) q_{j-1,\ell}(t) \, d\mu_{j-1}(t) =
q_{j-1,\ell}(r_{n-j+1}) \int q_{j,i}(t) \, d\mu_{j-1}(t).
$$
Thus,
$$
c_\ell \int q_{j-1,\ell}(t)^2 \, d\mu_{j-1}(t) = c_0
q_{j-1,\ell}(r_{n-j+1})  \int d\mu_{j-1}(t).
$$
Because $\ell \le i \le n-j$, both integrals are positive by
Lemma~\ref{lemma:posdefbd}. The largest root of $q_{j-1,\ell}$ is
less than $r_{n-j+1}$, so
$$
q_{j-1,\ell}(r_{n-j+1})>0.
$$
Thus, $c_0,\dots,c_i$ all have the same sign, which is positive
because $c_i=1$.
\end{proof}

Theorem~\ref{theorem:orthog} follows from applying
Lemma~\ref{lemma:induction} repeatedly, starting with $q_{k,n-k}$.

\section{Linear programming bounds}
\label{sec:lpbd}

Theorem~\ref{theorem:main} is deduced from the following
proposition, which originates in \cite{Y} and is the key technique
in \cite{Y,KY1,KY2,A1,A2} (here we give a somewhat different proof
from the one in \cite{Y}, along the lines of \cite{A2}):

\begin{proposition} \label{prop:lpbound}
Let $f \co (0,4] \to \R$ be any function.  Suppose $h \co [-1,1]
\to \R$ is a polynomial such that
$$
h(t) \le f(2-2t)
$$
for all $t \in [-1,1)$, and suppose there are nonnegative
coefficients $\alpha_0,\dots,\alpha_d$ such that $h$ has the
expansion
$$
h(t) = \sum_{i=0}^d \alpha_i C^{n/2-1}_i(t)
$$
in terms of ultraspherical polynomials.  Then every set of $N$
points on $S^{n-1}$ has potential energy at least
$$
N^2\alpha_0 - Nh(1)
$$
with respect to the potential function $f$.
\end{proposition}

Fundamentally, this proposition is a generalization of the linear
programming bounds for spherical codes due independently to
Kabatiansky and Levenshtein in \cite{KL} and Delsarte, Goethals,
and Seidel in \cite{DGS} (see also Chapter~9 in \cite{CS}).  As is
pointed out in \cite{Y}, one can derive those bounds for spherical
codes with minimal angle $\theta$ from
Proposition~\ref{prop:lpbound} by setting
$$
f(2-2t) = \begin{cases} \infty & \textup{if $t>\cos \theta$, and}\\
0 & \textup{otherwise.}
\end{cases}
$$
(Of course, $f$ takes values in $\R \cup \{\infty\}$, but that is
not a problem.)  The potential energy is then $0$ for a spherical
code with minimal angle at least $\theta$ and $\infty$ otherwise;
Proposition~\ref{prop:lpbound} implies that no such code can exist
if $N > h(1)/\alpha_0$.

Choosing the optimal function $h$ amounts to solving an
infinite-dimensional linear programming problem: the potential
energy bound is a linear functional of $h$, and the only
restrictions on $h$ are the linear constraints on its values and
ultraspherical coefficients. As in the case of the traditional
linear programming bounds, it is not hard to approximate the
optimum numerically (by solving a finite-dimensional linear
programming problem in which one imposes the constraint $h(t) \le
f(2-2t)$ for only finitely many values of $t$).  Even though the
bound is usually not sharp, it is often quite close.  For example,
for twenty points on $S^2$ and the Coulomb potential function
$f(r)=1/r^{1/2}$, the minimal-energy configuration, which
incidentally is not the vertices of a regular dodecahedron,
appears to have potential energy $301.763\dots$. The polynomial
\begin{eqnarray*}
h(t) &= & 0.8729 C^{1/2}_0(t) + 0.634 C^{1/2}_1(t) + 0.425
C^{1/2}_2(t) + 0.258 C^{1/2}_3(t)\\
&& \phantom{} + 0.135 C^{1/2}_4(t) + 0.0569 C^{1/2}_5(t) + 0.016
C^{1/2}_6(t)
\end{eqnarray*}
proves a lower bound of $301.204$.  It is not the best possible
choice of $h$, but we have found no polynomial that beats $301.5$.

What is remarkable is that in the case of the sharp arrangements,
for every completely monotonic potential function there is an
auxiliary function that proves a sharp bound, as we show in the
proof of Theorem~\ref{theorem:main}.

\begin{proof}[Proof of Proposition~\ref{prop:lpbound}]
Suppose $\CodeS \subset S^{n-1}$ has cardinality $N$.  By
assumption,
$$
\sum_{x,y \in \CodeS, x \ne y} f\big(|x-y|^2\big) \ge \sum_{x,y
\in \CodeS, x \ne y} h(\langle x,y \rangle),
$$
because $|x-y|^2 = 2-2\langle x,y \rangle$.  On the other hand,
\begin{eqnarray*}
\sum_{x,y \in \CodeS, x \ne y} h(\langle x,y \rangle) &=&
-Nh(1) + \sum_{x,y \in \CodeS} h(\langle x,y \rangle)\\
&=&  - Nh(1) + \sum_{i=0}^d \alpha_i \sum_{x,y \in \CodeS}
C^{n/2-1}_i(\langle x,y \rangle).
\end{eqnarray*}
Because $C^{n/2-1}_i$ is a positive-definite kernel (see
Subsection~\ref{subsec:pdk}),
\begin{equation*}
\sum_{x,y \in \CodeS} C^{n/2-1}_i(\langle x,y \rangle) \ge 0.
\displaybreak\end{equation*}
It follows that
$$
\sum_{x,y \in \CodeS, x \ne y} h(\langle x,y \rangle) \ge  - Nh(1)
+ \alpha_0 \sum_{x,y \in \CodeS} C^{n/2-1}_0(\langle x,y \rangle)
= N^2\alpha_0 - Nh(1),
$$
as desired.
\end{proof}

The bound proved by an auxiliary polynomial $h$ is sharp for a
configuration $\CodeS$ and potential function $f$ if and only if
two conditions hold: $h(t)$ must equal $f(2-2t)$ at every inner
product $t$ that occurs between distinct points in $\CodeS$, and
whenever the ultraspherical coefficient $\alpha_i$ is positive
with $i>0$, we must have
$$
\sum_{x,y \in \CodeS} C^{n/2-1}_i(\langle x,y \rangle) = 0.
$$
In particular, if $h$ is strictly positive definite, then $\CodeS$
must be a spherical $\deg(h)$-design, by Theorem~5.5 in
\cite{DGS}.

One can easily generalize Proposition~\ref{prop:lpbound} to deal
with the case of unequally charged particles (this was pointed out
in Theorem~1 of \cite{Y}). Suppose we have $N$ real numbers
$c_1,\dots,c_N$, which we refer to as charges; they must all have
the same sign.  Given a potential function $f$ and an ordered
$N$-tuple $(x_1,\dots,x_N) \in (S^{n-1})^N$, define the potential
energy to be
$$
\sum_{1 \le i,j \le N, i \ne j} c_ic_jf\big(|x_i-x_j|^2\big).
$$
Under the conditions in Proposition~\ref{prop:lpbound}, no such
configuration can have potential energy less than
$$
\left(\sum_{i=1}^N c_i\right)^2\alpha_0 - \left(\sum_{i=1}^N c_i^2
\right)h(1).
$$
We know of few beautiful configurations with unequal charges and
none for which linear programming bounds are sharp, but they may
well exist.

\section{Choosing auxiliary functions}
\label{sec:aux}

In this section, we describe an explicit choice for the auxiliary
function $h$ in Proposition~\ref{prop:lpbound};  it is analogous
to the constructions using the Christoffel-Darboux formula in
\cite{KL,MRRW}.  We will not require this function later in the
paper, but it provides a simple way to apply
Proposition~\ref{prop:lpbound}, and the techniques used in
analyzing it will be important in the proof of
Theorem~\ref{theorem:main}.

Let $f \co (0,4] \to \R$ be completely monotonic.  From the proof
of Proposition~\ref{prop:lpbound} we know that it is more natural
to look at $a(t) = f(2-2t)$, because $|x-y|^2 = 2-2\langle x,y
\rangle$.  The function $a$ is absolutely monotonic on $[-1,1)$
because $f$ is completely monotonic on $(0,4]$, and $a$ is
strictly absolutely monotonic if and only if $f$ is strictly
completely monotonic.

Our construction of $h$ requires three inputs: the dimension $n-1$
of the sphere $S^{n-1}$, a natural number $m \ge 1$, and a real
number $\alpha \ge 0$. Let
$$
t_1 < \dots < t_m
$$
be the roots of the polynomial $C_m^{n/2-1} + \alpha
C_{m-1}^{n/2-1}$.  We require that $t_1 \ge -1$ and thus
$\{t_1,\dots,t_m\} \subset [-1,1)$; by Theorem~3.3.4 in \cite{Sz},
that amounts to $\alpha \le -C_m^{n/2-1}(-1)/C_{m-1}^{n/2-1}(-1)$
(which is positive because $C^{n/2-1}_i(-1)$ has sign $(-1)^i$).

Let $h(t)$ be the Hermite interpolating polynomial that agrees
with $a(t)$ to order $2$ at each $t_i$ (i.e., $h(t_i)=a(t_i)$ and
$h'(t_i)=a'(t_i)$).  The idea of constructing auxiliary functions
via Hermite interpolation originates in \cite{Y}.

\begin{lemma} \label{lemma:upperbound1}
For all $t \in [-1,1)$,
$$
h(t) \le a(t).
$$
\end{lemma}

\begin{proof}
By Lemma~\ref{lemma:remainder}, there exists a point $\xi$ such
that
$$
\min(t,t_1,\dots,t_m) < \xi < \max(t,t_1,\dots,t_m)
$$
and
$$
a(t) - h(t) = \frac{a^{(2m)}(\xi)}{(2m)!}(t-t_1)^{2}\dots
(t-t_m)^{2}.
$$
The right side is nonnegative, and hence $a(t) \ge h(t)$.
\end{proof}

The last requirement before we can apply
Proposition~\ref{prop:lpbound}, namely that $h$ must be a
positive-definite function, is more subtle. Let
$$
F(t) = C_m^{n/2-1}(t)+ \alpha C_{m-1}^{n/2-1}(t) = \prod_{i=1}^m
(t-t_i).
$$
In the notation introduced before
Proposition~\ref{prop:absmonred}, $h = H(a,F^2)$.  We will show
that $F^2$ has a strong property that we refer to as conductivity:

\begin{definition}
A nonconstant polynomial $g$ with all its roots in $[-1,1)$ is
\textit{conductive\/} if for all absolutely monotonic functions
$a$ on $[-1,1)$, $H(a,g)$ is positive definite.  It is
\textit{strictly conductive\/} if it is conductive and for all
strictly absolutely monotonic $a$, $H(a,g)$ is strictly positive
definite of degree $\deg(g)-1$.
\end{definition}

\begin{lemma} \label{lemma:conductivity}
If $g_1$ and $g_2$ are conductive and $g_1$ is positive definite,
then $g_1g_2$ is conductive.  If $g_1$ is conductive and strictly
positive definite and $g_2$ is strictly conductive, then $g_1g_2$
is strictly conductive.
\end{lemma}

\begin{proof}
As noted in the proof of Proposition~\ref{prop:absmonred}, for
every absolutely monotonic function $a$ on $[-1,1)$,
$$
H(a,g_1g_2) = H(a,g_1) + g_1 H(Q(a,g_1),g_2).
$$
The lemma now follows from Proposition~\ref{prop:absmonred} and
the fact that the set of positive-definite (or strictly
positive-definite) functions is closed under taking products.
\end{proof}

For $r \in [-1,1)$, let $\ell_r$ denote the linear polynomial
$\ell_r(t) = t-r$.  Clearly $\ell_r$ is conductive, and even
strictly conductive if $r \ne -1$, because $H(a,\ell_r)$ is the
constant polynomial $a(r)$, which is nonnegative.
Lemma~\ref{lemma:conductivity} implies that if $g$ is conductive
and positive definite, then $g \ell_r$ is conductive.

Consider the partial products
$$
\prod_{i=1}^j (t-t_i)
$$
for $j \le m$.  By Theorem~\ref{theorem:orthog}, this product is
positive definite for $j<m$ (even strictly positive definite), and
for $j=m$ it is positive definite because $\alpha \ge 0$.  It
follows from Lemma~\ref{lemma:conductivity} and the remarks in the
previous paragraph that each product is conductive.  Thus,
$$
\prod_{i=1}^m (t-t_i)
$$
is conductive, so Lemma~\ref{lemma:conductivity} implies that its
square is as well.  Because $h=H(a,F^2)$, we conclude that $h$ is
positive definite, as desired.

These auxiliary polynomials are rarely optimal.  For example,
consider the case $n=3$, $N=20$, and $f(r)=1/r^{1/2}$ that was
discussed in Section~\ref{sec:lpbd}.  Numerical calculations
suggest that the best bound is obtained by taking $m=3$ and
$\alpha = 0.752718117\dots$ (an algebraic number of degree $22$).
That yields a lower bound of $299.708\dots$, which is worse than
the bound of $301.204$ derived in Section~\ref{sec:lpbd}. However,
it is reasonably close, especially given how simple and systematic
the construction of $h$ was.

The connection between our construction and those used in
\cite{MRRW} and \cite{KL} may not be immediately apparent. Those
constructions are equivalent to using the polynomial
$$
(t-t_m)\prod_{i=1}^{m-1}(t-t_i)^2
$$
in the linear programming bounds for spherical codes with minimal
angle $\theta$ satisfying $\cos \theta = t_m$.  The fact that
$$
\prod_{i=1}^{m-1}(t-t_i)
$$
is positive definite is proved in those papers using the
Christoffel-Darboux formula.  Unlike in our construction, $\alpha$
can be arbitrarily large.

One can generalize the construction from this section as follows.
Let $p_0,p_1,\dots$ be any family of orthogonal polynomials such
that each is positive definite.  Given $m$ and $\alpha$, choose
$t_1,\dots,t_m$ to be the roots of $p_m + \alpha p_{m-1}$.  The
construction and proof are essentially the same as before.  The
proof of Theorem~\ref{theorem:main} will be based on this
approach.

\section{Proof of the main theorem}
\label{sec:proof}

Let $f \co (0,4] \to \R$ be completely monotonic, define $a(t) =
f(2-2t)$ as in Section~\ref{sec:aux}, and let $\CodeS \subset
S^{n-1}$ be a sharp arrangement with $|\CodeS|=N$. To prove
Theorem~\ref{theorem:main}, we will construct an auxiliary
polynomial $h$ that satisfies the hypotheses of
Proposition~\ref{prop:lpbound} and proves a sharp bound.  The case
when $\CodeS$ is the set of vertices of the $600$-cell is
qualitatively more difficult than the other cases, so for the time
being we assume that we are in the sharp case.  (We deal with the
$600$-cell in Section~\ref{sec:600}.)

Suppose that the possible inner products that occur between
distinct points in $\CodeS$ are $t_1,\dots,t_m$, ordered so that
$$
-1 \le t_1 < \dots < t_m < 1.
$$
Let $h(t)$ be the Hermite interpolating polynomial that agrees
with $a(t)$ to order $2$ at each $t_i$.

When $t_1=-1$, it might seem more natural to interpolate only to
first order at that point.  The purpose of second-order
interpolation is to avoid sign changes in $h(t)-a(t)$ when
$t=t_i$, and that is not a concern when $t=-1$.  In fact, 
first-order interpolation works there, but we will omit the details
here, because the only advantage we see in it is that it leads to
a slightly improved bound on the number of derivative sign
conditions one must impose on the potential function (it decreases
by $1$ the bound $K$ mentioned in Section~\ref{sec:intro}). See
Section~\ref{sec:otherhomog} for the details of that approach in a
context in which it is more useful.

\begin{lemma} \label{lemma:upperbound}
For all $t \in [-1,1)$,
$$
h(t) \le a(t).
$$
\end{lemma}

The proof is identical to that of Lemma~\ref{lemma:upperbound1}.
Thus, the first hypothesis of Proposition~\ref{prop:lpbound} is
satisfied.  What remains to be shown is that $h$ is a
positive-definite function and that the bound is sharp.  (The
bound derived from this function $h$ is not sharp for the
$600$-cell, which is why that case is more subtle and we have
postponed it until later.)

The sharpness of the lower bound is simple, because it follows
from the definition of a sharp arrangement that $\CodeS$ is a
spherical $\deg(h)$-design.  Thus,
$$
\sum_{x,y \in \CodeS} C^{n/2-1}_i(\langle x,y \rangle) = 0
$$
whenever $0 < i \le \deg(h)$;  in fact, for each fixed $y$, the
polynomial $C^{n/2-1}_i(\langle x,y \rangle)$ has the same average
over $\CodeS$ as over the entire sphere $S^{n-1}$, and it
integrates to $0$ over the sphere because it is a positive-degree
spherical harmonic and thus orthogonal to the constant functions.
Because $h$ agrees with $a$ at the inner products $t_1,\dots,t_m$
from $\CodeS$, the conditions for a sharp bound that were listed
after the proof of Proposition~\ref{prop:lpbound} are satisfied.

The proof that $h$ is positive definite is analogous to the proof
in Section~\ref{sec:aux}, but it will involve more elaborate
lemmas.  Define
$$
F(t) = \prod_{i=1}^m (t-t_i).
$$
Then $h = H(a,F^2)$, and we will prove that $F^2$ is conductive.
In fact, we will show that $F^2$ is strictly conductive, which
will be important in the proof of uniqueness.

\begin{lemma} \label{lemma:gposdef}
The function $F$ is strictly positive definite.
\end{lemma}

This lemma is a standard fact in the theory of linear programming
bounds.  We prove it as in \cite{E}.

\begin{proof}
The leading coefficient of $F$ in terms of powers of $t$ is
positive, which implies that its leading ultraspherical
coefficient is also positive.  For the others, we apply
orthogonality to see that the $i$-th ultraspherical coefficient of
$F$ equals a positive constant (depending on $i$) times
$$
\int_{S^{n-1}} F(\langle x, y \rangle) C_i^{n/2-1}(\langle x, y
\rangle)\, dx,
$$
where $y$ is an arbitrary point on $S^{n-1}$ and of course $dx$
denotes surface measure on $S^{n-1}$.  Choose $y \in \CodeS$.

Because we know that the leading ultraspherical coefficient in $F$
is positive, we can take $i < \deg(F)$.  Thus, $\deg(F) + i \le
2\deg(F)-1$.  Because $\CodeS$ is a sharp configuration, it is a
spherical $(2\deg(F)-1)$-design.  It follows that the integral
above equals a positive constant times
$$
\sum_{x \in \CodeS} F(\langle x, y \rangle) C_i^{n/2-1}(\langle x,
y \rangle).
$$
By the construction of $F$, each term in the sum vanishes except
for $x=y$, so the sum equals $F(1) C_i^{n/2-1}(1)$.  Both factors
are positive, as desired.
\end{proof}

In order to deal with the partial products
$$
\prod_{i=1}^j (t-t_i),
$$
we must express $F$ in terms of orthogonal polynomials. Recall
that the ultraspherical polynomials $C^{n/2-1}_i$ are orthogonal
with respect to the measure
$$
(1-t^2)^{(n-3)/2}\, dt
$$ on $[-1,1]$.  Call that measure $d\mu(t)$.
Let $p_0,p_1,\dots$ denote the monic orthogonal polynomials with
respect to $(1-t) \,d\mu(t)$.  These polynomials are a special
case of Jacobi polynomials, but we will require only one fact
about them: they are nonnegative linear combinations of the
ultraspherical polynomials $C_i^{n/2-1}$.  That follows from
Lemma~\ref{lemma:signedmeasure} and
Proposition~\ref{prop:basecase}.

\begin{lemma} \label{lemma:orthoexpress}
There exists a constant $\alpha$ such that $F = p_m + \alpha
p_{m-1}$.
\end{lemma}

\begin{proof}
We simply need to show that $F$ is orthogonal to all polynomials
of degree at most $m-2$ with respect to the measure $(1-t)
\,d\mu(t)$.  That is equivalent to showing that $(1-t)F(t)$ is
orthogonal to all such polynomials with respect to $d\mu(t)$.

Let $p$ be any polynomial of degree at most $m-2$.  Because
$\CodeS$ is a sharp configuration, it must be a spherical
$(2m-1)$-design.  It follows as in the previous proof that
$$
\int (1-t)F(t) p(t) \, d\mu(t)
$$
equals a positive linear combination of
$$
(1-t_1)F(t_1) p(t_1),\dots,(1-t_m)F(t_m) p(t_m),(1-1)F(1) p(1).
$$
Each of these vanishes because $F$ vanishes at $t_1,\dots,t_m$.
Thus,
$$
\int F(t) p(t) \, (1-t)d\mu(t) = 0,
$$
as desired.
\end{proof}

Combining Lemma~\ref{lemma:orthoexpress} with
Theorem~\ref{theorem:orthog} shows that for $j < m$, the
polynomial
$$
\prod_{i=1}^j (t-t_i)
$$
is strictly positive definite, and the case $j=m$ is
Lemma~\ref{lemma:gposdef}. Now the proof that $F$ and $F^2$ are
conductive concludes as in Section~\ref{sec:aux}, by using
Lemma~\ref{lemma:conductivity} repeatedly.

We must still prove the additional uniqueness results when $f$ is
strictly completely monotonic.  For that, we use the following
lemma:

\begin{lemma}
If $a$ satisfies $a^{(k)}(t) > 0$ for all $k \ge 0$ and $t \in
(-1,1)$, then $a(t)-h(t)$ has at most $\deg(h)+1$ roots in
$[-1,1)$, counted with multiplicity.
\end{lemma}

\begin{proof}
By Rolle's theorem, $a(t)-h(t)$ has at most one more root than
$a'(t)-h'(t)$.  Repeated differentiation reduces to the case in
which $\deg(h)=0$, which is trivial.
\end{proof}

Thus, the only roots of $a(t)-h(t)$ in $[-1,1)$ are at $t=t_i$ for
some $i$. It follows from the proof of
Proposition~\ref{prop:lpbound} that if $\CodeT$ is any other
potential energy minimum with $|\CodeT|=|\CodeS|$, then all inner
products between distinct points of $\CodeT$ occur among
$t_1,\dots,t_m$. Thus $\CodeT$ is a spherical code with the same
parameters as $\CodeS$, and each sharp configuration listed in
Table~\ref{tab:sharp} except on the last line is known to be the
unique spherical code with its parameters, up to orthogonal
transformations, of course. See Appendix~\ref{app} for details and
references.

Even when the spherical code is not unique, one can still conclude
from the strict absolute monotonicity of $a$ that $\CodeT$ must be
sharp.  The key observation is that $F^2$ is strictly conductive
(except in the trivial case when $\CodeS$ consists of two
antipodal points), which follows easily from
Lemma~\ref{lemma:conductivity}. It implies that $h$ is a strictly
positive-definite polynomial of degree $2m-1$. Then the sharpness
in Proposition~\ref{prop:lpbound} implies that $\CodeT$ must be a
spherical $\deg(h)$-design (by Theorem~5.5 in \cite{DGS}) and
hence a sharp configuration.  Furthermore, each of $t_1,\dots,t_m$
must occur as an inner product in $\CodeT$, from switching the
roles of $\CodeS$ and $\CodeT$.

\section{The $600$-cell}
\label{sec:600}

The final configuration is the vertices of the regular $600$-cell.
The construction in Section~\ref{sec:proof} does not work in this
case, and overall the $600$-cell appears to be intrinsically more
complicated than the sharp configurations.  The fundamental
problem is that it is only a spherical $11$-design, but the
polynomial $h$ constructed as in the previous section would have
degree $15$ (or $14$ if one uses the alternate construction for
antipodal configurations). Recall that being a spherical
$\deg(h)$-design was crucial to having a sharp bound. In fact, $h$
does turn out to be a positive-definite function in this case, but
it proves a suboptimal bound.

Andreev and Elkies have shown in \cite{A3,E} that a degree $17$
polynomial with coefficients in $\Q(\sqrt{5})$ proves a sharp
bound for the number of points in a spherical code with minimal
angle $\pi/5 = \cos^{-1} (1+\sqrt{5})/4$ on $S^3$, although
polynomials of lower degree do not work. A similar phenomenon
occurs in our problem. What saves the proof is that although the
$600$-cell is not a spherical $12$-design, all spherical harmonics
of degrees from $13$ to $19$ do indeed sum to $0$ over the
$600$-cell.  Degree $12$ is the only problem, and that can be
addressed by ensuring that the $12$-th ultraspherical coefficient
of $h$ vanishes. The fact that spherical harmonics of degrees $13$
through $19$ sum to $0$ over the $600$-cell is easily checked
using the distance distribution from Table~\ref{tab:dist}, which
tells for each vertex how many others have a given inner product
with it. One can also prove it using invariant theory as follows.
The symmetry group of the $600$-cell is the $H_4$ reflection
group, for which the ring of invariant polynomials has generators
of degrees $2$ (the trivial invariant $x \mapsto |x|^2$), $12$,
$20$, and $30$. Summing a polynomial over the vertices of the
$600$-cell amounts to averaging it with respect to $H_4$.  For any
degree that is not a nonnegative integer combination of $12$,
$20$, and $30$, averaging must yield a polynomial that is constant
on $S^3$, and that constant vanishes for spherical harmonics
because they are orthogonal to the constant polynomials.

\begin{table}
\caption{The distance distribution of the $600$-cell.}\label{tab:dist}
\begin{center}
\begin{tabular}{cc}
Inner Product & Count\\
\hline $\pm 1$ & $1$\\
$(\pm 1 \pm \sqrt{5})/4$ & $12$\\
$\pm 1/2$ & $20$\\
$0$ & $30$\\
\\ 
\end{tabular}
\end{center}
\end{table}

For the $600$-cell, we have $m=8$ and $\{t_1,\dots,t_m\} = \{-1,
0, \pm 1/2, (\pm 1 \pm \sqrt{5})/4\}$; we order the inner products
so that $t_1 < \dots < t_8$. Let $h(t)$ be the unique polynomial
of degree at most $17$ such that $h(t_i) = a(t_i)$ for $1 \le i
\le 8$, $h'(t_i) = a'(t_i)$ for $2 \le i \le 8$, and $\alpha_{11}
= \alpha_{12} = \alpha_{13} = 0$, where $\alpha_i$ denotes the
$i$-th ultraspherical coefficient of $h$. Note that we do not
require that $h'(-1) = a'(-1)$.

It might seem more natural to use a polynomial of degree $15$ and
simply set $\alpha_{12}=0$, but that does not work.  We have no
conceptual explanation for why it fails or why our slightly more
complicated approach succeeds.

If $h(t) \le a(t)$ for all $t$ and $h$ is positive definite, then
it proves a sharp bound (as pointed out above, because
$\alpha_{12}=0$). Unfortunately, we know of no simple reason why
these inequalities hold.  We have proved them using a computer
algebra system.  The proof is completely rigorous, in the sense
that it uses no floating point arithmetic or other unrigorous
approximations;  all computations are performed using exact
arithmetic in $\Q(\sqrt{5})$.

\subsection{Proof that $h$ is positive definite}

We wish to prove that $h(t)$ is a nonnegative linear combination
of the ultraspherical polynomials $C^1_i(t)$.  Each ultraspherical
coefficient of $h(t)$ is a linear function of $a(t_i)$ and
$a'(t_i)$ for $1 \le i \le m$.  Let $u_i$ denote the coefficient
of $a(t_i)$ in this linear function, and let $v_i$ denote the
coefficient of $a'(t_i)$.

As explained in Section~\ref{sec:intro}, it suffices to consider
the potential functions $f(r) = (4-r)^k$ with $k \in
\{0,1,2,\dots\}$; up to scaling, that amounts to taking $a(t) =
(1+t)^k$.  The linear combination
$$
\sum_{i=1}^m \big(u_i a(t_i) + v_i a'(t_i)\big)
$$
then becomes
$$
\sum_{i=1}^m \big(u_i (1+t_i)^k + v_i k (1+t_i)^{k-1}\big).
$$
If $v_8>0$, then this expression is positive for all sufficiently
large $k$.  To prove that it is nonnegative for all $k$, we simply
compute a bound $\ell$ such that it is guaranteed to be positive
for all $k \ge \ell$ and then check nonnegativity for $0 \le k <
\ell$.

Let
$$
\chi(x) = \begin{cases} x & \textup{if $x\le 0$, and}\\
0 & \textup{if $x>0$.}
\end{cases}
$$
If $\ell$ is chosen so that
\begin{equation}
\label{eq:ellsum} \sum_{i=1}^m \big(\chi(u_i) (1+t_i)^\ell +
\chi(v_i) \ell (1+t_i)^{\ell-1}\big) + v_8 \ell (1+t_8)^{\ell-1}
\ge 0,
\end{equation}
then
$$
\sum_{i=1}^m \big(u_i (1+t_i)^k + v_i k (1+t_i)^{k-1}\big) \ge 0
$$
for all $k \ge \ell$.  In other words, once the asymptotically
dominant term outweighs all the negative terms, it continues to do
so forever (because whenever $k$ increases, the dominant term
increases by a larger factor than any other term does).

To prove that $h$ is positive definite, we deal with eighteen
cases, one for each ultraspherical coefficient of $h$.  In each
case, we compute the coefficients
$$
u_1,\dots,u_8,v_1,\dots,v_8
$$
explicitly as elements of $\Q(\sqrt{5})$.  We then check that
$v_8>0$ and that \eqref{eq:ellsum} holds with $\ell=32$.  Finally,
we check that the ultraspherical coefficient is nonnegative for $k
\in \{0,1,\dots,31\}$.  These calculations are too cumbersome to
carry out by hand, but they are easily done using a computer
algebra system.

\subsection{Proof that $h(t) \le a(t)$}

Unfortunately, the inequality $h(t) \le a(t)$ no longer has as
simple a proof as in the case of the sharp configurations.  To
prove it, we will use Proposition~\ref{prop:absmonred}.  Let
$$
F(t) = (t+1)\prod_{i=2}^m (t+t_i)^2,
$$
and let $\tilde{h}$ be the usual Hermite interpolation $H(a,F)$ of
$a$ (without requiring any ultraspherical coefficients to vanish).
By Proposition~\ref{prop:absmonred},
$$
\frac{a(t) - \tilde{h}(t)}{F(t)}
$$
is absolutely monotonic on $[-1,1)$.  Let $q(t)$ be the quadratic
Taylor polynomial for this quotient about $t=-1$.  Then it follows
from absolute monotonicity that
$$
\frac{a(t) - \tilde{h}(t)}{F(t)} \ge q(t)
$$
for all $t \in [-1,1)$.  Thus,
$$
\frac{a(t) - h(t)}{F(t)} \ge q(t) +
\frac{\tilde{h}(t)-h(t)}{F(t)}.
$$
The right side of this inequality is a quadratic polynomial.  If
we verify that it is nonnegative over $[-1,1)$, then $a(t) \ge
h(t)$ over that range (because $F(t) \ge 0$).  To check this, we
simply check that it is nonnegative at $t= \pm 1$ and has
nonpositive leading coefficient. Each of these calculations
amounts to checking that a certain explicit linear combination of
$a(t_1),\dots,a(t_8)$, $a'(t_1),\dots,a'(t_8)$, $a^{(2)}(-1)$, and
$a^{(3)}(-1)$ is nonnegative.  We prove it using the method from
the previous subsection, with $\ell=36$.

\subsection{Proof of uniqueness}

Finally, we must deal with the issue of uniqueness when the
potential function $f$ is strictly completely monotonic (so $a$ is
strictly absolutely monotonic).  In that case, by
Proposition~\ref{prop:absmonred},
$$
\frac{a(t) - \tilde{h}(t)}{F(t)}
$$
is strictly absolutely monotonic on $(-1,1)$, which implies that
$$
\frac{a(t) - \tilde{h}(t)}{F(t)} > q(t)
$$
for $t \in (-1,1)$.  Thus, $a(t)-h(t)$ has roots only at $t=t_i$
for some $i$. Uniqueness follows, because the vertices of the
regular $600$-cell are the only $(4,120,(1+\sqrt{5})/4)$ spherical
code (see \cite[p.~260]{B} or \cite{BD}).

\section{Other compact two-point homogeneous spaces}
\label{sec:otherhomog}

The key property underlying our proofs is the two-point
homogeneity of the sphere. Recall that a metric space is
\textit{two-point homogeneous\/} if for each $r>0$, its isometry
group acts transitively on ordered pairs of points at distance
$r$.  In this section we generalize Theorem~\ref{theorem:main} to
all compact, connected two-point homogeneous spaces.  This
generalization is parallel to the general setting for linear
programming bounds in \cite{CS,L2,L3}. It seems plausible that our
results also generalize to the case of discrete two-point
homogeneous spaces (see Chapter~9 of \cite{CS} for examples), but
we have not investigated that possibility.

Let $X$ be a compact, connected two-point homogeneous space with
more than one point. Wang proved in \cite{Wa} that, up to changing
to an equivalent metric (i.e., applying an invertible function to
the metric), $X$ is one of the following possibilities: a sphere,
a real projective space, a complex projective space, a
quaternionic projective space, or the octonionic projective plane,
each with its standard metric.  In other words, we use the round
metric on the sphere, the Fubini-Study metric on real, complex, or
quaternionic projective space, and the $F_4$-invariant metric on
$\Oc\Proj^2$ (up to scaling, they are the unique two-point
homogeneous Riemannian metrics on these spaces).

The octonionic projective plane $\Oc\Proj^2$ may be the least
familiar of these possibilities; it cannot be defined as
$\Oc^3\setminus\{0\}$ modulo the action of $\Oc^*$ because of the
nonassociativity of the octonions ($\Oc^*$ is not a group and
does not act on $\Oc^3$). However, it can be defined using the
exceptional Jordan algebra or as $F_4/\textup{Spin}(9)$. There is
no such space as $\Oc\Proj^3$, because $\Oc\Proj^2$ does not
satisfy the Desargues theorem.  See \cite{Ba} for more details on
octonionic projective geometry.

Normalize the metric $d$ on $X$ (the geodesic distance with
respect to the Riemannian metric) so that the greatest distance
between points is $\pi$. Just as we use inner products to measure
distances between points on spheres, it is natural to use $\cos
d(x,y)$ instead of $d(x,y)$.  In $\R\Proj^n$, $\C\Proj^n$, or
$\Ha\Proj^n$ this function is easy to calculate: if we represent
$x$ and $y$ by unit vectors in $\R^{n+1}$, $\C^{n+1}$, or
$\Ha^{n+1}$, respectively, then $\cos d(x,y) = 2|\langle x,y
\rangle|^2-1$, where the brackets denote the standard Hermitian
form over $\R$, $\C$, or $\Ha$.

Spheres are $n$-point homogeneous for every $n$, but the other
two-point homogeneous spaces are not even three-point homogeneous.
One consequence is that configurations in projective space can
frequently be deformed without changing the distances between
points.  For example, for $\tau \in \C$ with
$|\tau|=(1+\sqrt{5})/2$, consider the configuration
$\mathcal{C}_\tau$ consisting of the six points in $\C\Proj^2$
with homogeneous coordinates $[\tau,1,0]$, $[\tau,-1,0]$,
$[1,0,\tau]$, $[-1,0,\tau]$, $[0,\tau,1]$, and $[0,\tau,-1]$. When
$\tau = (1+\sqrt{5})/2$, they are the images in $\R\Proj^2 \subset
\C\Proj^2 $ of the vertices of a regular icosahedron in $\R^3$.
The cosine of the distance between any two distinct points is
$-3/5$, so the points form a regular simplex embedded in
$\C\Proj^2$.  It will follow from Theorem~\ref{theorem:projective}
that each $\mathcal{C}_\tau$ is a universally optimal six-point
code in $\C\Proj^2$. Lemma~\ref{lemma:equiv} shows that these
configurations form a one-dimensional family.  By contrast,
positive-dimensional families of sharp configurations do not exist
in spheres.

\begin{lemma} \label{lemma:equiv}
For $|\tau|=|\tau'|=(1+\sqrt{5})/2$, there is an isometry of
$\C\Proj^2$ mapping $C_\tau$ to $C_{\tau'}$ if and only if
$\tau/\tau'$ or $\tau/\bar\tau'$ is a sixth root of unity.
\end{lemma}

\begin{proof}
The isometry group of $\C\Proj^2$ is generated by the action of
$U(3)$ on lines in $\C^3$ together with complex conjugation.
Complex conjugation sends $C_\tau$ to $C_{\bar\tau}$, the diagonal
unitary transformation that multiplies the $j$-th coordinate by
$e^{2\pi i j/3}$ sends $C_\tau$ to $C_{\tau e^{-2\pi i/3}}$, and
the identity map sends $C_\tau$ to $C_{-\tau}$.  In the other
direction, if $x$, $y$, and $z$ are unit vectors in the six lines
in $\C^3$ corresponding to the points of $\mathcal{C}_\tau$, then
$$
\langle x,y \rangle \langle y,z \rangle \langle z,x \rangle \in
\left\{1,1/5,\pm 1/\big(5\sqrt{5}\big),
\pm{\tau^3}/{\big(|\tau|^2+1\big)^3},
\pm{\bar\tau^3}/{\big(|\tau|^2+1\big)^3}\right\},
$$
and each of these values occurs (note that the triple product
$\langle x,y \rangle \langle y,z \rangle \langle z,x \rangle$ is
independent of which unit vector is chosen in each line). Because
$|\tau|^3/{\big(|\tau|^2+1\big)^3} = 1/\big(5\sqrt{5}\big)$, the
isometry class of $\mathcal{C}_\tau$ determines the set $\{\pm 1,
\pm\tau^3,\pm \bar\tau^3\}$, and that proves the necessity of the
condition in the lemma statement.
\end{proof}

The potential functions we will consider are the absolutely
monotonic functions of $\cos d(x,y)$.  As in the case of the
sphere, these functions can also be understood in terms of another
metric, namely the \textit{chordal distance\/} $d_c(x,y) =
\sin\big(d(x,y)/2\big)$ introduced for Grassmannians in \cite{CHS}
(i.e., $\sqrt{1-|\langle x,y\rangle|^2}$ when $x$ and $y$ are
represented by unit vectors as above). For motivation, recall that
the Euclidean distance between two points separated by an angle of
$\theta$ on the unit sphere is $2\sin (\theta/2)$. For projective
spaces, the chordal distance is most naturally expressed in terms
of the orthogonal projection operators $\Pi_x$ and $\Pi_y$ onto
the lines $x$ and $y$.  Then
$$
d_c(x,y) = \frac{|\Pi_x - \Pi_y|_F}{\sqrt{2}},
$$
where $|\cdot|_F$ is the Frobenius norm, defined by
$$
|A|_F^2 = \sum_{i,j} |A_{i,j}|^2 = \tr \bar{A}^t A.
$$
Absolutely monotonic functions of $\cos d(x,y)$ are the same as
completely monotonic functions of squared chordal distance, just
as in the case of the unit sphere.

Let $G$ be the isometry group of $X$.  Then one can decompose
$L^2(X)$ as
$$
L^2(X) = \widehat{\bigoplus_{\ell \ge 0}} \,\, V_\ell,
$$
where $V_0,V_1,\dots$ are certain finite-dimensional irreducible
representations of $G$. This decomposition is unique because the
summands are pairwise nonisomorphic. For each $\ell$, one can
define a positive-definite kernel $K_\ell$ on $X$ as in the case
of the sphere, and one can define a 
function $C_\ell \co [-1,1] \to \R$ such that
$$
K_\ell(x,y) = C_\ell(\cos d(x,y)).
$$
We index the representations $V_\ell$ so that $C_\ell$ is a
polynomial of degree $\ell$.  One can compute these polynomials
explicitly (see for example \cite[p.~178]{G}). For projective
spaces, the result is that $C_\ell$ is a positive constant times
the Jacobi polynomial $P^{(\alpha,\beta)}_\ell$, where $\alpha =
n/2-1$, $n$ is the dimension of $X$ as a real manifold, and
$\beta$ is $-1/2$, $0$, $1$, or $3$ according to whether one is
working over $\R$, $\C$, $\Ha$, or $\Oc$, respectively (i.e.,
$\beta = (\dim_\R A)/2-1$ for the algebra $A$). We define a
\textit{positive-definite} polynomial to be a nonnegative linear
combination of $C_0,C_1,\dots$.

An \textit{$M$-design\/} in $X$ is a finite subset $\CodeS \subset
X$ such that for $f \in V_i$ with $1 \le i \le M$,
$$
\sum_{x \in \CodeS} f(x) = 0.
$$
As in the case of spheres, being an $M$-design is also
characterized by the requirement that
$$
\sum_{x,y \in \CodeS} C_i(\cos d(x,y)) = 0
$$
for $1 \le i \le M$.  For example, to verify that the
configurations $\mathcal{C}_\tau$ constructed above are
$1$-designs in $\C\Proj^2$, one need only check that $6 C_1(1) +
30 C_1(-3/5)=0$.

We say that a finite subset $\CodeS \subset X$ is a \textit{sharp
configuration\/} if it is an $M$-design, $m$ distances occur
between distinct points in $\CodeS$, and $M \ge 2m-1-\delta$,
where $\delta=1$ if there are two antipodal points in $\CodeS$ and
$\delta=0$ otherwise.  By a pair of ``antipodal points'' we mean
points separated by the greatest possible distance in $X$; this
condition is more subtle in projective spaces than in spheres,
because it is no longer true that each point has a unique
antipode.

This definition of a sharp configuration agrees with the previous
one in the case of spheres. Here we merely require an antipodal
sharp configuration to be a $(2m-2)$-design, but in the case of
spheres every antipodal $(2m-2)$-design is automatically a
$(2m-1)$-design: by Theorem~7.4 of \cite{DGS}, every point in such
a design has an antipode if one does, and that suffices to make
every odd-degree homogeneous polynomial sum to zero.  On the other
hand, some antipodal sharp configurations in projective space are
only $(2m-2)$-designs.  For example, the $D_4$, $E_6$, and $E_7$
root systems give antipodal sharp configurations of $12$, $36$,
and $63$ points in $\R\Proj^3$, $\R\Proj^5$, and $\R\Proj^6$,
respectively.  For these codes $m=2$ (the cosines of the distances
are $-1$ and $-1/2$), and they are $2$-designs but not
$3$-designs. By contrast, $E_8$ yields a $3$-design in
$\R\Proj^7$.

The following theorem generalizes Theorem~\ref{theorem:main} and
the results of Levenshtein from \cite{L2} (namely, that all sharp
configurations are optimal codes):

\begin{theorem} \label{theorem:projective}
Let $f \co [-1,1) \to \R$ be absolutely monotonic, and let $\CodeS
\subset X$ be a sharp configuration. If $\CodeT \subset X$ is any
subset satisfying $|\CodeT|=|\CodeS|$, then
$$
\sum_{x,y \in \CodeT, x \ne y} f(\cos d(x,y)) \ge \sum_{x,y \in
\CodeS, x \ne y} f(\cos d(x,y)).
$$
If $f$ is strictly absolutely monotonic, then equality implies
that $\CodeT$ is a sharp configuration and that the same distances
occur between points in $\CodeS$ and $\CodeT$.
\end{theorem}

Naturally, we call
$$
\sum_{x,y \in \CodeS, x \ne y} f(\cos d(x,y))
$$
the $f$-potential energy of $\CodeS$.  Strictly speaking, this
definition conflicts with our earlier definition for the special
case $X=S^{n-1}$, because previously we used functions of squared
Euclidean distance as our potential functions.  However, that
should not cause any confusion: the Euclidean distance is twice
the chordal distance, so this change simply amounts to the
difference between using the squared chordal distance and the
cosine of the geodesic distance. In the terminology of
Section~\ref{sec:proof}, here we are dealing with the absolutely
monotonic function $a$ instead of the completely monotonic
potential function, and thus our previous framework is equivalent
to the current one in the case of $S^{n-1}$.

The proof of Theorem~\ref{theorem:projective} is parallel to the
proof for spheres, so we merely sketch it here, with particular
attention to arguments that require modification in projective
space.

The linear programming bounds generalize in a straightforward manner, with
essentially the same proof:

\begin{proposition} \label{prop:projlpbound}
Let $f \co [-1,1) \to \R$ be any function.  Suppose $h \co [-1,1]
\to \R$ is a polynomial such that $h(t) \le f(t)$ for all $t \in
[-1,1)$, and suppose there are nonnegative coefficients
$\alpha_0,\dots,\alpha_d$ such that $h$ has the expansion
$$
h(t) = \sum_{i=0}^d \alpha_i C_i(t).
$$
Then every set of $N$ points in $X$ has $f$-potential energy at
least
$$
N^2\alpha_0 - Nh(1).
$$
\end{proposition}

We apply this proposition in the same way as before. Suppose first
that $\CodeS$ is a $(2m-1)$-design, with $m$ distances between
distinct points in $\CodeS$.  Let $t_1 < \dots < t_m$ be the
cosines of those distances.  As before, we construct the Hermite
interpolating polynomial $h$ of degree $2m-1$ that agrees with $f$
to order $2$ at each of $t_1,\dots,t_m$.  The proof that $h(t) \le
f(t)$ for all $t \in [-1,1)$ is the same as in
Lemma~\ref{lemma:upperbound}.

Define
$$
F(t) = \prod_{i=1}^m (t-t_i).
$$
The same proof as in Lemma~\ref{lemma:gposdef} shows that $F$ is
positive definite, and we apply an analogue of
Lemma~\ref{lemma:orthoexpress} to express $g$ in terms of a
related family of orthogonal polynomials as $p_m+\alpha p_{m-1}$.
To complete the proof that $h$ is positive definite, we use the
fact that the polynomials in that family are themselves positive
definite, together with Theorem~\ref{theorem:orthog}, to establish
conductivity. Everything is completely parallel to the spherical
case. Finally, because $\CodeS$ is a $(2m-1)$-design and
$\deg(h)=2m-1$, the lower bound for potential energy must be
sharp.  When $f$ is strictly absolutely monotonic, the uniqueness
argument is the same as before.

When $t_1=-1$ and $\CodeS$ is a $(2m-2)$-design but not a
$(2m-1)$-design, one must instead use first-order Hermite
interpolation at $t=-1$, which complicates the proof. We can
assume that $m>1$, since if $m=1$, then every pair of points in
$\CodeS$ is antipodal and Theorem~\ref{theorem:projective} visibly
holds.

If we define
$$
F(t) = \prod_{i=2}^m (t-t_i),
$$
then the heart of the proof is the verification that $(t+1)F(t)^2$
is strictly conductive, and everything else works exactly the same
as before.  Let $q_0,q_1,\dots$ be the monic orthogonal
polynomials with respect to $(1-t^2)\,d\mu(t)$, where $\mu$ is the
measure on $[-1,1]$ with respect to which the polynomials $C_\ell$
are orthogonal.  (Recall that they are certain Jacobi polynomials
$P^{(\alpha,\beta)}_\ell$.  Then $d\mu(t) =
(1-t)^\alpha(1+t)^\beta \, dt$.)  The argument used to prove
Lemma~\ref{lemma:orthoexpress} shows that $F = q_{m-1} + \alpha
q_{m-2}$ for some $\alpha$.

The most subtle part of the proof of
Theorem~\ref{theorem:projective} is the proof that $\alpha > 0$.
(If $\CodeS$ were a $(2m-1)$-design, then $\alpha$ would vanish.)
We start with a crucial observation about the polynomials $q_i$:

\begin{lemma} \label{lemma:qposdef}
For each $i$, the polynomial $(t+1)q_i(t)$ is strictly positive
definite.
\end{lemma}

\begin{proof}
This polynomial is orthogonal to all polynomials of degree less
than $i$ with respect to $(1-t) \, d\mu(t)$, so for some constant
$c$ it equals $p_{i+1}(t) + c p_i(t)$ in terms of the monic
orthogonal polynomials for this measure.  Setting $t=-1$ shows
that $c \ge 0$, because each $p_j$ has all its roots in $(-1,1)$
and hence $p_j(-1)$ has sign $(-1)^j$. As before, it follows from
Lemma~\ref{lemma:signedmeasure} and
Proposition~\ref{prop:basecase} that each $p_j$ is strictly
positive definite.
\end{proof}

\begin{lemma} \label{lemma:ineqs}
Let $f(t)$ be a monic polynomial, with $0$-th coefficient $f_0$ in
terms of $C_0(t),C_1(t),\dots$.
\begin{enumerate}
\item If $\deg(f) \le 2m-2$ and $f(t_i)=0$ for $1 \le i \le m$,
then $f_0|\CodeS| = f(1)$. \item If $\deg(f) = 2m-1$ and
$f(t_i)=0$ for $1 \le i \le m$, then $f_0|\CodeS| < f(1)$.\item If
$\deg(f) \le 2m-2$ and $f(t) \ge 0$ for $-1 \le t \le 1$, then
$f_0|\CodeS| \ge f(1)$.
\end{enumerate}
\end{lemma}

\begin{proof}
Write
$$
f(t) = \sum_{i=0}^{\deg(f)} f_i C_i(t).
$$
In cases (1) and (2),
$$
\sum_{x,y \in \CodeS} f\big(\cos d(x,y) \big) = |\CodeS| f(1),
$$
while in case (3),
$$
\sum_{x,y \in \CodeS} f\big(\cos d(x,y) \big) \ge |\CodeS| f(1).
$$
On the other hand,
$$
\sum_{x,y \in \CodeS} f\big(\cos d(x,y) \big) =
\sum_{i=0}^{\deg(f)} f_i \sum_{x,y \in \CodeS} C_i\big(\cos
d(x,y)\big).
$$
The $i=0$ term yields $f_0 |\CodeS|^2$, while all terms with $1
\le i \le 2m-2$ vanish because $\CodeS$ is a $(2m-2)$-design.  In
case (1) that yields $f_0|\CodeS| = f(1)$ and in case (3) it
yields $f_0|\CodeS| \ge f(1)$.  Finally, in case (2) we have
$f_{2m-1}>0$ because $f$ is monic, so the $i=2m-1$ term is
positive (because $\CodeS$ is not a $(2m-1)$-design) and hence
$f_0|\CodeS| < f(1)$.
\end{proof}

\begin{lemma} \label{lemma:zerocoeff}
If $f$ is a polynomial with $\deg(f) \le i$, then
$(t+1)f(t)q_i(t)$ has $0$-th coefficient equal to $f(1)$ times
that of $(t+1)q_i(t)$ (in terms of $C_0(t),C_1(t),\dots$).
\end{lemma}

\begin{proof}
The $0$-th coefficient of a polynomial $p(t)$ equals
$$
\frac{\int p(t) \, d\mu(t)}{\int d\mu(t)}.
$$
Divide $f(t)$ by $t-1$ to yield $f(t) = g(t)(t-1) + f(1)$, with
$\deg(g) < i$.  Then
$$
\int (t+1)f(t)q_i(t) \, d\mu(t) = \int g(t) q_i(t) \,
(1-t^2)d\mu(t) + f(1) \int (1+t)q_i(t) \, d\mu(t).
$$
Because $q_i$ is orthogonal to all polynomials of degree less than
$i$ with respect to $(1-t^2)d\mu(t)$,
$$
\int (t+1)f(t)q_i(t) \, d\mu(t) = f(1) \int (1+t)q_i(t) \,
d\mu(t),
$$
which completes the proof.
\end{proof}

Let $\beta_i$ denote the $0$-th coefficient of $(t+1)q_{i}(t)$ in
terms of $C_0(t),C_1(t),\dots$. Taking $f(t) = (t+1)F(t)$ in case
(1) of Lemma~\ref{lemma:ineqs}, we find that
$$
(\beta_{m-1} + \beta_{m-2}\alpha)|\CodeS| = 2q_{m-1}(1) +
2q_{m-2}(1)\alpha.
$$
For an upper bound on $|\CodeS|$, take $f(t) =
(t+1)F(t)q_{m-1}(t)$ in case (2). By Lemma~\ref{lemma:zerocoeff},
$$
|\CodeS| < \frac{f(1)}{F(1)\beta_{m-1}} =
\frac{2q_{m-1}(1)}{\beta_{m-1}}
$$
(because $F(1)\beta_{m-1}>0$). For a lower bound, take $f(t) =
(t+1)q_{m-2}(t)^2$ in case (3). By Lemma~\ref{lemma:zerocoeff},
$$
|\CodeS| \ge \frac{f(1)}{q_{m-2}(1)\beta_{m-2}} =
\frac{2q_{m-2}(1)}{\beta_{m-2}}
$$
(because $q_{m-2}(1)\beta_{m-2}>0$).  We have shown that
$$
2q_{m-1}(1) - \beta_{m-1}|\CodeS| > 0
$$
and
$$
2q_{m-2}(1) - \beta_{m-2}|\CodeS| \le 0,
$$
while
$$
2q_{m-1}(1) - \beta_{m-1}|\CodeS| + \alpha \big(2q_{m-2}(1) -
\beta_{m-2}|\CodeS|\big) = 0.
$$
It follows that $\alpha>0$.

To prove that $(t+1)F(t)^2$ is strictly conductive using
Lemma~\ref{lemma:conductivity}, we need only show that all its
partial products
$$
(t+1)\prod_{i=2}^j (t-t_i)^2
$$
and
$$
(t+1)(t-t_{j+1})\prod_{i=2}^j (t-t_i)^2
$$
are strictly positive definite. That follows from
Theorem~\ref{theorem:orthog}, Lemma~\ref{lemma:qposdef}, and the
fact that the set of nonnegative linear combinations of
$q_0,q_1,q_2,\dots$ is closed under taking products (see
\cite{Gas} for the most general case or Theorem~5.2 in
\cite[p.~44]{A} for a simpler proof that applies to these
orthogonal polynomials).  Once the conductivity of $(t+1)F(t)^2$
has been established, the rest of the proof of
Theorem~\ref{theorem:projective} is essentially identical to that
of Theorem~\ref{theorem:main}.

\section{The Euclidean case} \label{sec:Euclidean}

One can generalize Proposition~\ref{prop:lpbound} to apply to
configurations in Euclidean space.  Here we outline this approach,
which leads to an optimization problem reminiscent of the
Beurling-Selberg extremal functions (see \cite{V}). We deal only
with the case of periodic, discrete point configurations (unions
of finitely many translates of a lattice), to avoid worries about
how potential energy is defined in pathological cases.

Define the \textit{density\/} of a periodic arrangement to be the
number of points per unit volume in space.  If the arrangement
consists of the $N$ translates $\Lambda+v_1,\dots,\Lambda+v_N$ of
a lattice $\Lambda \in \R^n$, with $v_j-v_k \not\in \Lambda$ for
$j \ne k$, then the density is $N/\vol(\R^n/\Lambda)$.  Given a
potential function $f \co (0,\infty) \to [0,\infty)$, define the
\textit{potential energy\/} of the packing to be
$$
\frac{1}{N}\sum_{1 \le j,k \le N} \sum_{x \in \Lambda, x+v_j-v_k
\ne 0} f\big(|x+v_j-v_k|^2\big).
$$
Of course summing over all pairs of points in the packing would
give a divergent sum (even this sum may diverge, in which case we
say the potential energy is infinite). The potential energy is the
average over all points in the configuration of the sum over all
distances to other points. In particular, the same periodic
arrangement can be written in terms of different $\Lambda$ and
$v_1,\dots,v_N$, but the potential energy is independent of this
choice.

For a different point of view, imagine placing finitely many
points $v_1,\dots,v_N$ on the flat torus $\R^n/\Lambda$. The
energy of a point is given by a sum over all geodesics connecting
it to the other points or to itself, and the energy of the
configuration is the average of this quantity over all the points.

Let $B_R(0)$ denote the ball of radius $R$ about $0$. Given a
configuration $\mathcal{C} \subset \R^n$, let
$$
\mathcal{C}_R = \mathcal{C} \cap B_R(0).
$$

\begin{lemma} \label{lemma:limit}
Let $\mathcal{C}$ be a periodic point configuration in $\R^n$. If
the potential energy of $\mathcal{C}$ is finite, then it equals
$$
\lim_{R \to \infty} \frac{1}{|\mathcal{C}_R|} \sum_{x,y \in
\mathcal{C}_R, x \ne y} f\big(|x-y|^2\big).
$$
The same holds even without the assumption that $f$ is
nonnegative, as long as $f(x) = O\big((1+|x|)^{-n/2-\delta}\big)$
for some $\delta>0$.
\end{lemma}

\begin{proof}
Suppose $\mathcal{C}$ is the disjoint union of
$\Lambda+v_1,\dots,\Lambda+v_N$ as above.  Say a point in
$\mathcal{C}$ has type $j$ if it is in $\Lambda+v_j$.  With each
point $x\in \mathcal{C}$ we associate the sum
\begin{equation}\label{eq:partial}
\sum_{y \in \mathcal{C}, y \ne x} f\big(|x-y|^2\big).
\end{equation}
This sum depends only on the type of $x$, and the definition of
the $f$-potential energy of $\mathcal{C}$ equals the average over
all types of this quantity.

All these sums must converge, because of the bounds on $f$ if $f$
is not assumed nonnegative or else because of the hypothesis that
the potential energy is finite. Let $\varepsilon>0$ be small, and
choose $K$ such that for $x$ of each type,
$$
\sum_{y \in \mathcal{C}, |y-x| > K} \big|f\big(|x-y|^2\big)\big|
\le \varepsilon.
$$
In particular, summing in \eqref{eq:partial} only over $y$ such
that $|x-y| \le K$ yields a partial sum within $\varepsilon$ of
the total sum.

In the sum
$$
\sum_{x,y \in \mathcal{C}_R, x \ne y} f\big(|x-y|^2\big),
$$
we can focus our attention on the points $x \in B_{R-K}(0)$, since
the number in $B_R(0) \setminus B_{R-K}(0)$ is negligible compared
to the number in $B_{R-K}(0)$ (on the order of $R^{n-1}$ rather
than $R^n$ as $R \to \infty$).  Each point contributes a bounded
amount to the sum, so omitting the points in $B_R(0) \setminus
B_{R-K}(0)$ changes the sum by $O\big(R^{n-1}\big)$.

For each $x \in \mathcal{C}_{R-K}$,
$$
\sum_{y \in \mathcal{C}_R, y \ne x} f\big(|x-y|^2\big)
$$
is within $\varepsilon$ of
$$
\sum_{y \in \mathcal{C}, y \ne x} f\big(|x-y|^2\big).
$$
The numbers of points of each type in $\mathcal{C}_{R-K}$ are
equal, up to a factor of $1+O(1/R)$.  Thus, for $R$ large,
$$
\frac{1}{|\mathcal{C}_R|} \sum_{x,y \in \mathcal{C}_R, x \ne y}
f\big(|x-y|^2\big)
$$
is $O(1/R)$ plus a quantity within $\varepsilon$ of the
$f$-potential energy of $\mathcal{C}$. Choosing $\varepsilon$
arbitrarily small and $R$ accordingly large completes the proof.
\end{proof}

The Euclidean analogue of Proposition~\ref{prop:lpbound} is stated
in terms of the Fourier transform on $\R^n$. We normalize the
Fourier transform of an $L^1$ function $h \co \R^n \to \R$ by
$$
\widehat{h}(t) = \int_{\R^n} h(x) e^{2\pi i \langle x,t \rangle}
\, dx.
$$

\begin{lemma} \label{lemma:integral}
Let $\mathcal{C}$ be a periodic point configuration in $\R^n$,
with density $\delta$. If $\varepsilon>0$ is sufficiently small
(depending on $\mathcal{C}$), then
$$
\lim_{R \to \infty} \int_{B_\varepsilon(0)} \frac{\left|\sum_{x
\in \mathcal{C}_R} e^{2\pi i \langle x,t \rangle}
\right|^2}{|\mathcal{C}_R|} \, dt = \delta.
$$
\end{lemma}

\begin{proof}
We will prove that for $\varepsilon>0$ sufficiently small, if $g$
is any radial, smooth function supported in $B_\varepsilon(0)$,
then
$$
\lim_{R \to \infty} \int_{\R^n} \frac{\left|\sum_{x \in
\mathcal{C}_R} e^{2\pi i \langle x,t \rangle}
\right|^2}{|\mathcal{C}_R|} g(t) \, dt = \delta g(0).
$$
Because the characteristic function of a small ball can be bounded
above and below by such functions, that will suffice to prove the
lemma.

We expand the numerator of the integrand via
$$
\left|\sum_{x \in \mathcal{C}_R} e^{2\pi i \langle x,t \rangle}
\right|^2 = \sum_{x,y \in \mathcal{C}_R} e^{2\pi i \langle x-y,t
\rangle}.
$$
It follows that
$$
\int_{\R^n} \frac{\left|\sum_{x \in \mathcal{C}_R} e^{2\pi i
\langle x,t \rangle} \right|^2}{|\mathcal{C}_R|} g(t) \, dt =
\frac{1}{|\mathcal{C}_R|}\sum_{x,y \in \mathcal{C}_R}
\widehat{g}(x-y).
$$
Let $\mathcal{C}$ be the disjoint union of translates
$\Lambda+v_1,\dots,\Lambda+v_N$ of a lattice $\Lambda$.  It
follows from Lemma~\ref{lemma:limit} that
$$
\lim_{R \to \infty} \frac{1}{|\mathcal{C}_R|}\sum_{x,y \in
\mathcal{C}_R} \widehat{g}(x-y) = \frac{1}{N} \sum_{1 \le j,k \le
N} \sum_{z \in \Lambda} \widehat{g}(z+v_j-v_k).
$$
(This quantity is not quite the same as that considered in
Lemma~\ref{lemma:limit}, because we allow $x=y$ or $z+v_j-v_k=0$
in the sum, but the extra terms amount to $\widehat{g}(0)$ on each
side.)

Recall the Poisson summation formula, which states that for every
Schwartz function $h \co \R^n \to \R$ and every $v \in \R^n$,
$$
\sum_{z \in \Lambda} h(z+v) = \frac{1}{\vol(\R^n/\Lambda)} \sum_{t
\in \Lambda^*} \widehat{h}(t) e^{-2\pi i \langle v,t \rangle},
$$
where
$$
\Lambda^* = \{ t \in \R^n : \langle t,z \rangle \in \Z \textup{
for all $z \in \Lambda$} \}
$$
is the dual lattice.  (See Appendix~A of \cite{CK2} for the
standard proof.) Taking $h = \widehat{g}$ yields
$$
\frac{1}{N} \sum_{1 \le j,k \le N} \sum_{z \in \Lambda}
\widehat{g}(z+v_j-v_k) = \frac{1}{N\vol(\R^n/\Lambda)} \sum_{t \in
\Lambda^*} g(t) \left|\sum_{1 \le j \le N} e^{2\pi i \langle v_j,
t\rangle }\right|^2.
$$
If the support of $g$ is sufficiently small, then only the $t=0$
term contributes to the right side, and it equals $\delta g(0)$
because $N/\vol(\R^n/\Lambda) = \delta$.
\end{proof}

\begin{proposition} \label{prop:Euclidean}
Let $f \co (0,\infty) \to [0,\infty)$ be any function. Suppose $h
\co \R^n \to \R$ satisfies $h(x) \le f\big(|x|^2\big)$ for all $x
\in \R^n \setminus \{0\}$ and is the Fourier transform of a
function $g \in L^1(\R^n)$ such that $g(t) \ge 0$ for all $t \in
\R^n$. Then every periodic configuration in $\R^n$ with density
$\delta$ has $f$-potential energy at least
$$
\delta \big(\liminf_{t \to 0} g(t)\big)-h(0).
$$
\end{proposition}

Without loss of generality, one can assume that $g$ and $h$ are
radial functions (replace $g$ with the average of its rotations
about the origin). For simplicity of notation, we typically write
$g = \widehat{h}$. That is justified because $g$ is in fact the
Fourier transform of $h$ as a tempered distribution (when these
functions are radial; of course in general $g(t) =
\widehat{h}(-t)$), but note that we have not assumed that $h$ is
integrable and thus it is not clear that
$$
\int_{\R^n} h(x) e^{2\pi i \langle x,t \rangle} \, dx
$$
converges.  By contrast, $h(x)$ does equal
$$
\int_{\R^n} \widehat{h}(t) e^{2\pi i \langle x,t \rangle} \, dt.
$$

\begin{proof}
Define
$$
\ell = \liminf_{t \to 0} g(t).
$$
Let $\varepsilon>0$, and choose $\eta>0$ such that $g(t) \ge
\ell-\varepsilon$ whenever $|t| \le \eta$.

Let $\mathcal{C}$ be a periodic configuration with density
$\delta$.  By Lemma~\ref{lemma:limit}, the $f$-potential energy is
$$
\lim_{R \to \infty} \frac{1}{|\mathcal{C}_R|} \sum_{x,y \in
\mathcal{C}_R, x \ne y} f\big(|x-y|^2\big).
$$
Because $h(x) \le f\big(|x|^2\big)$ for all $x \ne 0$,
$$
\frac{1}{|\mathcal{C}_R|} \sum_{x,y \in \mathcal{C}_R, x \ne y}
f\big(|x-y|^2\big) \ge -h(0) + \frac{1}{|\mathcal{C}_R|} \sum_{x,y
\in \mathcal{C}_R} h(x-y).
$$
On the other hand,
\begin{eqnarray*}
\frac{1}{|\mathcal{C}_R|} \sum_{x,y \in \mathcal{C}_R} h(x-y) &=&
\frac{1}{|\mathcal{C}_R|} \sum_{x,y \in \mathcal{C}_R} \int_{\R^n}
g(t) e^{2\pi i \langle x-y, t \rangle} \, dt\\
&=& \int_{\R^n} g(t) \frac{\left|\sum_{x \in \mathcal{C}_R}
e^{2\pi i \langle x,t
\rangle} \right|^2}{|\mathcal{C}_R|} \, dt\\
&\ge& \int_{B_\eta(0)} (\ell-\varepsilon) \frac{\left|\sum_{x \in
\mathcal{C}_R} e^{2\pi i \langle x,t \rangle}
\right|^2}{|\mathcal{C}_R|} \, dt.
\end{eqnarray*}
In the limit as $R \to \infty$, Lemma~\ref{lemma:integral} implies
that the potential energy is at least
$$
-h(0) + \delta(\ell - \varepsilon),
$$
if $\eta$ is sufficiently small.  Because $\varepsilon$ can be
taken to be arbitrarily small, the desired lower bound follows.
\end{proof}

This proposition is closely analogous to the linear programming
bounds for sphere packing introduced in \cite{CE} (see also
\cite{C,CK1,CK2}).  The proof given here is somewhat different,
because we have arranged it to minimize the hypotheses on $h$.  If
one is willing to assume stronger hypotheses, namely bounds on the
decay rates of $h$ and $\widehat{h}$, then one can give a proof
closer to that in \cite{CE} by moving the Poisson summation from
the proof of Lemma~\ref{lemma:integral} into that of
Proposition~\ref{prop:Euclidean}.

Analogously to the case of Proposition~\ref{prop:lpbound}, one can
derive the bounds from \cite{CE} (with weaker hypotheses than in
that paper) by setting
$$
f(r) = \begin{cases} \infty & \textup{for $r < 1$, and}\\
0 & \textup{otherwise.}
\end{cases}
$$
Then one concludes that the sphere centers in a packing with balls
of radius $1/2$ must have density at most $h(0)/\widehat{h}(0)$
(assuming $\widehat{h}$ is continuous at $0$).

Let $\Lambda_2$, $\Lambda_8$, and $\Lambda_{24}$ denote the
hexagonal lattice in $\R^2$, the $E_8$ root lattice in $\R^8$, and
the Leech lattice in $\R^{24}$, respectively.

\begin{conjecture}\label{conj:opt}
Let $n \in \{2,8,24\}$, and let $f \co (0,\infty) \to \R$ be
completely monotonic and satisfy $f(x) =
O\big(|x|^{-n/2-\varepsilon}\big)$ as $|x| \to \infty$ for some
$\varepsilon>0$. Then there exists a function $h$ that satisfies
the hypotheses of Proposition~\ref{prop:Euclidean} and proves that
$\Lambda_n$ has the least $f$-potential energy of any periodic
configuration in $\R^n$ with its density.
\end{conjecture}

One consequence of this conjecture is that these lattices are
universally optimal configurations in Euclidean space.  That
appears to be a rare property, which provably fails for lattices
in dimensions $3$, $5$, $6$, and $7$. The unique optimal lattices
in those dimensions are the $A_3$, $D_5$, $E_6$, and $E_7$ root
lattices (see \cite{Mar}). Each of these lattices has higher
energy than its dual lattice for the potential function $f(x) =
e^{-x}$, when both lattices are rescaled to have unit density.

Montgomery proved in \cite{Mo} that the hexagonal lattice is
universally optimal among all lattices in $\R^2$, but that is a
weaker assertion than universal optimality among all periodic
point configurations.  Sarnak and Str\"ombergsson conjectured in
\cite{SS} that $D_4$, $E_8$, and the Leech lattice are universally
optimal among lattices and proved that they are local optima.
Perhaps one could combine their theorem with the techniques in
this paper and \cite{CK1} to prove that $E_8$ and $\Lambda_{24}$
are the global minima among lattices in their dimensions.

If true, Conjecture~\ref{conj:opt} implies that $E_8$ and the
Leech lattice are the densest sphere packings in their dimensions,
by using the potential function $f(r) = 1/r^s$ with $s$ large.
(The density of a sphere packing is the fraction of space covered
by the spheres; unfortunately this terminology conflicts with the
density of a point configuration.  We will refer to these numbers
as the ``packing density'' and ``point density'' to avoid
confusion.) The potential energy of any periodic packing is
asymptotic to
$$
\frac{K}{r_{\min}^{2s}}
$$
as $s \to \infty$, where $r_{\min}$ is the smallest distance
between distinct points and $K$ is the average number of points at
distance $r_{\min}$ from a point in the packing (i.e., the average
kissing number if one centers balls of radius $r_{\min}/2$ at the
points of the arrangement to form a sphere packing). Universal
optimality implies that $r_{\min}$ must be as large as possible
for a periodic configuration of point density $1$, which means the
packing density is maximized.  Thus, Conjecture~\ref{conj:opt}
implies that no periodic packing can exceed the packing density of
$E_8$ or the Leech lattice, and therefore no packing can (see
Appendix~A of \cite{CE}).

Furthermore, Conjecture~\ref{conj:opt} implies that $E_8$ and the
Leech lattice are the unique densest periodic packings in $\R^8$
and $\R^{24}$.  In the notation used above, any equally dense
periodic packing must have at least as large an average kissing
number $K$, or the potential energy would be lower.  Because $E_8$
and $\Lambda_{24}$ have the greatest possible kissing numbers in
their dimensions and those kissing configurations are unique (see
Chapters~13 and~14 of \cite{CS}), the kissing configurations in
the hypothetical packing are determined up to orthogonal
transformations.  However, those orthogonal transformations may be
nontrivial, so it is not obvious that the packing is uniquely
determined.

\begin{lemma} \label{lemma:unique}
For $n=8$ or $n=24$, there is a unique sphere packing in $\R^n$
such that each kissing configuration is the $\Lambda_n$ kissing
configuration.
\end{lemma}

\begin{proof}
Without loss of generality, suppose one sphere is centered at a
point in $\Lambda_n$ and the spheres tangent to it are centered
at the nearest points in $\Lambda_n$. We must prove that each
sphere is centered at a point of $\Lambda_n$. Call a sphere center
\textit{good\/} if it is in $\Lambda_n$ and the same is true for
all spheres tangent to it, and call one point a neighbor of
another if they are the centers of tangent spheres. By assumption,
at least one point is good; call it the origin.  We will prove
that all neighbors of good points are also good.  It follows that
all points in the packing are good.  To see why, suppose $x$ is
not good with $|x|$ minimal.  Then $x$ must have a neighbor $y$
with $|y| < |x|$, since otherwise a suitable scalar multiple of
$x$ (chosen between $x$ and the origin to match the minimal
distance of the packing) could be included with the neighbors of
$x$ to increase the kissing number in $\R^n$; the reason is that
every neighbor of $x$ would be on or outside the sphere of radius
$|x|$ and hence at least as far from the new point as $x$ is. The
smaller neighbor $y$ is good by assumption, which implies that $x$
is good.

It remains to prove that all neighbors of good points are good.
Suppose $x$ is good and $y$ is a neighbor of $x$.  We must show
that there is a unique way to arrange the neighbors of $y$
(namely, as points of $\Lambda_n$).  The positions of the common
neighbors of $x$ and $y$ are already determined, because $x$ is
good.  It would suffice to show that every isometry of $\R^n$ that
fixes $x$ and $y$ and preserves the set of common neighbors must
preserve $\Lambda_n$, because any two possible kissing
arrangements about $y$ are related by such an isometry.  (They are
related by some isometry fixing $y$ by the uniqueness of the
kissing configuration.  Because the symmetry group acts
transitively on the neighbors (see Lemma~1.10 in \cite[p.~25]{Eb}
and Theorem~27 in \cite[p.~288]{CS}), we may assume that $x$ is
fixed as well, and it follows that the set of common neighbors is
also preserved.) For $n=8$ we will show that this is true; for
$n=24$ it is false, but we will prove a slightly weaker statement
that will still suffice.

Without loss of generality, we can translate the origin so that
$y=0$. When $n=8$, $x$ and its $56$ common neighbors with $0$ span
$E_8$, so every orthogonal transformation that preserves $x$ and
the common neighbors clearly preserves $E_8$.  That completes the
proof by showing that the kissing arrangements about two neighbors
can only be in one orientation relative to each other.

When $n=24$, $x$ and its $4600$ common neighbors with $0$ do not
span $\Lambda_{24}$ (their span cannot contain any vector that has
an odd inner product with $x$).  Instead, they span a sublattice
$L$ of index~$2$, as is not hard to check using Figure~4.12 in
\cite[p.~133]{CS} (with $x = (4,-4,0,\dots,0)/\sqrt{8}$). This
lattice contains all the minimal vectors from $\Lambda_{24}$ that
have even inner product with $x$, and it contains twice the
minimal vectors from $\Lambda_{24}$ that have inner product $\pm
1$ with $x$.

The dual lattice
$$
L^* = \{ w \in \R^n : \langle w,z \rangle \in \Z \textup{ for all
$z \in L$} \}
$$
contains $L$ (because $L$ is integral), and considering
covolumes shows that $[L^* : L] = 4$.  The Leech lattice is
contained in $L^*$, as is the vector $x/2$, and they must generate
$L^*$ because they generate a lattice in which $L$ has index $4$.

Every orthogonal transformation that preserves $x$ and the set of
$4600$ common neighbors preserves $L$ and hence $L^*$.  That means
the kissing configuration of our packing at $0$ must be contained
among the vectors of norm $4$ in $L^*$, because there exists such
an orthogonal transformation that yields the Leech kissing
configuration, which is contained in $L^*$. Besides the minimal
vectors of $\Lambda_{24}$, the only other possibilities are the
vectors of the form $\pm(x/2 - z)$ where $z$ is a minimal vector
of $\Lambda_{24}$ satisfying $\langle x,z \rangle = 1$. (Any
vector in $L^*$ but not $\Lambda_{24}$ must be of the form $x/2-w$
with $w \in \Lambda_{24}$.  A short calculation shows that if
$|x/2-w|^2=4$, then either $w=z$ with $|z|^2=4$ and $\langle x,z
\rangle =1$ or $w=x-z$ for such a $z$.)

The kissing configuration must contain all the minimal vectors of
$L$, namely the minimal vectors of $\Lambda_{24}$ that have even
inner product with $x$, because $L$ is preserved under the
transformation to the Leech kissing configuration. It must then
either contain $\pm z$ for all $z \in \Lambda_{24}$ with $|z|^2=4$
and $\langle x,z \rangle = 1$, in which case it is simply all the
minimal vectors in $\Lambda_{24}$, or it must contain all the
vectors $\pm(x/2 - z)$ for such $z$. (If it contains one of the
points $\pm z$, then it contains all of them because any one of
them together with $L$ spans $\Lambda_{24}$ and hence the
orthogonal transformation would be a symmetry of $\Lambda_{24}$.)
One can check as follows that the second possibility is indeed a
valid kissing configuration about $0$. If $r$ denotes reflection
in the hyperplane through the origin orthogonal to $x$, then the
map $z \mapsto -r(z)$ preserves $L$, maps $x$ to $x$, and
interchanges the two kissing configurations.

We have shown that there are only two distinct kissing
configurations about $0$ that both contain $x$ and the $4600$
common neighbors.  Fortunately, only the first can occur in the
sphere packing.  Suppose the kissing configuration about $0$
contained $x/2-z$ with $z \in \Lambda_{24}$, $|z|^2=4$, and
$\langle x,z \rangle = 1$. Because $x$ is good, its kissing
configuration contains the point $x-z$, which is so close to
$x/2-z$ that the spheres centered at those points would overlap.
That contradiction rules out the second case, so $0$ must be good
as well, which completes the proof.
\end{proof}

Under Conjecture~\ref{conj:opt}, let $h_s$ prove a sharp bound for
$f(r)=1/r^s$.  It seems likely that, when suitably rescaled, $h_s$
converges as $s \to \infty$ to a function that proves a sharp
bound in the setting of \cite{CE}.  (One can check that this holds
for the functions constructed below when $n=1$.)  That would put
the mysterious functions from that paper into a broader context
and exhibit them as limits of families of functions.

We have not been able to prove Conjecture~\ref{conj:opt}, but we
can prove the one-dimensional analogue.  Of course, minimizing
potential energy in $\R^1$ is completely trivial, so this result
has no interesting geometric consequences.  However, it is not
obvious that Proposition~\ref{prop:Euclidean} proves a sharp bound
even in this trivial case, so we view its sharpness as an argument
in favor of Conjecture~\ref{conj:opt}.

\begin{proposition}\label{prop:onedim}
Let $f \co (0,\infty) \to \R$ be completely monotonic and satisfy
$f(x) = O\big(|x|^{-1/2-\varepsilon}\big)$ for some
$\varepsilon>0$ as $|x| \to \infty$. Then there exists a function
$h \co \R \to \R$ such that $h(x) \le f\big(x^2\big)$ for all $x
\in \R \setminus \{0\}$, $h$ is the Fourier transform of a
continuous function $\widehat{h} \co \R \to \R$ with support in
$[-1,1]$, $\widehat{h}(t) \ge 0$ for all $t \in \R$, and
$$
\widehat{h}(0)-h(0) = \sum_{k \in \Z\setminus \{0\}}
f\big(k^2\big).
$$
\end{proposition}

For notational convenience, define $\ft(x) = f\big(x^2\big)$ for
$x \in \R\setminus \{0\}$ (and let $\ft'(x) = 2xf'(x^2)$ denote
its derivative). Then we can take
\begin{eqnarray*}
h(x) &=& \sum_{k \in \Z\setminus\{0\}} \ft(k) \left(\frac{\sin
\pi(x-k)}{\pi (x-k)}\right)^2 + \sum_{k \in \Z\setminus\{0\}}
\ft'(k) \frac{\big(\sin \pi(x-k)\big)^2}{\pi^2(x-k)}\\
&& \phantom{} - 2 \left(\frac{\sin \pi x}{\pi x}\right)^2
\sum_{k=1}^\infty \big(\ft(k)+k\ft'(k)\big).
\end{eqnarray*}
(Of course we define $h(x)$ by continuity when $x \in \Z$.) By the
mean value theorem and the fact that $f'$ is increasing and
nonpositive,
$$
\left|\frac{f(x)-f(x/2)}{x/2}\right| \ge |f'(x)|,
$$
which implies that $|xf'(x)| = O\big(|x|^{-1/2-\varepsilon}\big)$;
thus, the sum of $k \ft'(k)$ in the definition of $h$ converges,
because $\ft'(k) = 2k f'\big(k^2\big)$.

We have no reason to think that these functions are the only ones
that prove a sharp bound in Proposition~\ref{prop:Euclidean} for
$n=1$, but they are the only ones we know how to construct.   We
prove below that they satisfy the properties listed in
Proposition~\ref{prop:onedim}.

The construction of $h$ is based on a variant of the Shannon
sampling theorem for band-limited functions (see Theorem~9 in
\cite{V} for the particular interpolation formula we use).  It is
natural to search for a function $h$ such that $\widehat{h}$
vanishes outside $[-1,1]$, by analogy with the one-dimensional
construction in Section~5 of \cite{CE}.  One can reconstruct such
a function from its values and derivatives at the integers,
assuming $h$ is integrable, via
$$
h(x) = \sum_{k \in \Z} h(k) \left(\frac{\sin \pi(x-k)}{\pi
(x-k)}\right)^2 + \sum_{k \in \Z} h'(k) \frac{\big(\sin
\pi(x-k)\big)^2}{\pi^2(x-k)}.
$$
We will not require this general fact, but it is important
motivation for the definition of $h$.

It is easy to check from the definition of $h$ that $h(k) =
\ft(k)$ and $h'(k) = \ft'(k)$ for $k \in \Z \setminus \{0\}$.
These conditions are necessary if $h$ is to prove a sharp bound in
Proposition~\ref{prop:Euclidean}.  We also have $h'(0)=0$ (in
fact, $h$ is even).  However, the value of $h(0)$, which is given
by
$$
h(0) = - 2 \sum_{k=1}^\infty \big(\ft(k) + k \ft'(k)\big),
$$
is more mysterious.  Setting $h(0)=\ft(0)$ would make no sense,
because $\ft(0)=f(0)$, which is frequently undefined.  The value
we have chosen for $h(0)$ makes $h$ have the right properties, but
we know no motivation for it that does not rely on some
calculation.

Because
$$
\left(\frac{\sin \pi x}{\pi x}\right)^2 = \int_{-1}^1 (1-|t|)
e^{2\pi i t x}\, dt
$$
and
$$
\frac{(\sin \pi x)^2}{\pi^2 x} = \int_{-1}^1 \frac{\sgn t}{2\pi i}
e^{2\pi i t x} \, dt,
$$
$h$ is the Fourier transform of a function $\widehat{h}$ supported
on $[-1,1]$.  Because $h$ is even, $\widehat{h}$ is also even, so
it suffices to compute $\widehat{h}(t)$ for $t \in [0,1]$.  (Of
course the simplification is inessential, but it avoids absolute
value signs.)  For such $t$, one finds that
\begin{eqnarray*}
\widehat{h}(t) &=& 2(1-t) \sum_{k=1}^\infty \ft(k) \cos 2\pi k t
- \frac{1}{\pi} \sum_{k=1}^\infty \ft'(k) \sin 2 \pi k t\\
&& \phantom{} - 2 (1-t) \sum_{k=1}^\infty \big(\ft(k) + k
\ft'(k)\big).
\end{eqnarray*}
One motivation for the value of $h(0)$ is that if one formally
computes $\widehat{h}'(1)$, one gets the answer $0$.  However, we
do not require this fact and make no attempt to justify it.

To complete the proof of Proposition~\ref{prop:onedim}, we must
prove that $h(x) \le f\big(x^2\big)$ for all $x \ne 0$ and
$\widehat{h}(t) \ge 0$ for all $t$.  One potential approach is to
reduce to the case of $f(x) = e^{-cx}$ with $c>0$: by Bernstein's
theorem (Theorem~12b in \cite[p.~161]{W}), there exists a weakly
increasing function $\alpha \co [0,\infty) \to \R$ such that for
all $x>0$,
$$
f(x) = \int_0^\infty e^{-cx} \, d\alpha(c).
$$
We can thus assume without loss of generality that $f(x) =
e^{-cx}$ with $c>0$.  (There is no contribution from $c=0$ because
if there were, then $f$ would not decay at infinity.) To see that
the inequalities we need are preserved under this reduction, we
must check that integration over $c$ commutes with the sums over
$k$ in the formulas for $h$ and $\widehat{h}$.  That follows from
Fubini's theorem.

We have not fully carried out this approach, but it is likely
possible.  For the less general case of $f(x) = e^{-c\sqrt{x}}$,
which corresponds to the special case when not only $f$ but also
$x \mapsto f\big(x^2\big)$ is completely monotonic, the inequality
$h(x) \le f\big(x^2\big)$ is very nearly proved in Theorem~8 of
\cite{GV}, which constructs similar functions to solve a related
optimization problem.  That paper does not contain precisely the
inequality we need, but a simple modification of the techniques
used there suffices, and one can furthermore prove that
$\widehat{h}(t) \ge 0$ by computing $\widehat{h}$ explicitly and
manipulating it. This proof suffices for potential functions $f$
such that $x \mapsto f\big(x^2\big)$ is completely monotonic, but
it is somewhat complicated and unilluminating.

Here we check these inequalities (and hence prove
Proposition~\ref{prop:onedim}) via a more conceptual approach, by
proving that they follow from the corresponding assertions about
codes on $S^1$ by taking limits.  Consider $N$ points on $S^1$,
under the completely monotonic potential function $g$;  we will
define $g$ in terms of $f$ shortly. In the proof of
Theorem~\ref{theorem:main}, we constructed an auxiliary function
$h_{N,g}$ that proved a sharp bound for the $g$-potential energy
in Proposition~\ref{prop:lpbound}.  When $N$ is odd, one can write
it down explicitly as follows.  Define $F(t) = (T_N(t)-1)/(t-1)$,
where $T_N$ is the $N$-th Chebyshev polynomial, given by $T_N(\cos
\theta) = \cos N\theta$.  (The roots of $F$ are the inner products
other than $1$ in the regular $N$-gon.)

\begin{lemma} \label{lemma:hformula}
When $N$ is odd,
$$
h_{N,g}(t) = \frac{F(t)}{N^2} \sum_{r : F(r)=0} \frac{a_r +
b_r(t-r)}{(t-r)^2},
$$
where
$$
a_r = 2(1-r)^2(1+r)g(2-2r)
$$
and
$$
b_r = -4(1-r)^2(1+r)g'(2-2r) - 2(1-r)(1+2r)g(2-2r).
$$
\end{lemma}

The sum is over all roots $r$ of $F$, each included only once
despite the fact that they are double roots. The proof of this
formula amounts to verifying that the right side is the Hermite
interpolation of $t \mapsto g(2-2t)$ to second order at the roots
of $F$, which is a straightforward calculation. (The case of even
$N$ differs slightly because one must treat the single root $r=-1$
of $F$ differently.  We will assume without loss of generality
that $N$ is odd, but the even case can easily be handled by the
same techniques.)

We define $g$ by
$$
g(x) = f\big(N^2 x/\big(4\pi^2\big)\big) .
$$
The motivation is that as $N \to \infty$,
$$
g\big(2 - 2\cos (2\pi x/N)\big) = \ft\big((N/\pi) \sin (\pi x /N)
\big) \to \ft(x),
$$
which suggests the \pagebreak following lemma.

\begin{lemma} \label{lemma:converge}
For $g(x) = f\big(N^2 x/\big(4\pi^2\big)\big)$, we have
$$
\lim_{N \to \infty} h_{N,g}\big(\cos (2\pi x/N)\big) = h(x)
$$
as $N \to \infty$ with $N$ odd.
\end{lemma}

The same is also true for $N$ even, but we will not need it. Note
that Lemma~\ref{lemma:converge} makes no assertion as to the
uniformity of the convergence (pointwise convergence will suffice
for our purposes).

\begin{proof}
At a formal level, it follows from taking limits term by term and
using
$$
\lim_{N \to \infty} \frac{g'\big(2 - 2\cos (2\pi
x/N)\big)}{\big(N/(2\pi)\big)^2} = \frac{\ft'(x)}{2x}
$$
and
$$
\lim_{N \to \infty} \frac{F\big(\cos(2\pi x/N)\big)}{N^2} =
\left(\frac{\sin \pi x}{\pi x} \right)^2
$$
that $h_{N,g}\big(\cos (2\pi x/N)\big)$ converges to
$$
\left(\frac{\sin \pi x}{\pi x} \right)^2 \sum_{k = 1}^\infty
\left(\ft(k)\left(\frac{2kx-k^2}{(x-k)^2}-\frac{2kx+k^2}{(x+k)^2}\right)
+ \ft'(k)\left(\frac{k^2}{x-k}-\frac{k^2}{x+k}\right) \right),
$$
which is equal to $h(x)$.  (The roots of $F$ are of the form
$\cos(2\pi k/N)$, which is how the sum over $r$ becomes a sum over
$k$.) To prove rigorously that
$$
\lim_{N \to \infty} h_{N,g}\big(\cos (2\pi x/N)\big) = h(x),
$$
one must justify taking the limit term by term.

To do so, we will apply the dominated convergence theorem.  Set
$$
r = \cos(2\pi k/N) = 1-2\sin^2(\pi k/N)
$$
for $0 < k < N/2$, and turn the sum over $r$ in the formula for
$h_{N,g}(t)$ into a sum over all $k>0$ by including zero terms for
all larger $k$. Let $t = 1-2\sin^2(\pi x/N)$.  We will prove that
the $k$-th term in the sum is $O\big(k^{-1-2\varepsilon}\big)$ in
absolute value, where the implicit constant is independent of $N$.
(Of course this calculation uses $f(x) =
O\big(|x|^{-1/2-\varepsilon}\big)$ and $f'(x) =
O\big(|x|^{-3/2-\varepsilon}\big)$.)  Because
$$
\sum_{k=1}^\infty k^{-1-2\varepsilon}
$$
converges, taking limits term by term is justified by the
dominated convergence theorem.

One can calculate that the $k$-th term is given by
\begin{equation}\label{eq:bigformula}
\frac{a_r + b_r(t-r)}{(t-r)^2} =
\frac{2\big(2s_k^2-1+3(s_x/s_k)^2-4s_x^2\big)}{\big(1-(s_x/s_k)^2\big)^2}
f_k -\frac{4N^2
s_k^2\big(1-s_k^2\big)}{\pi^2\big(1-(s_x/s_k)^2\big)}f'_k,
\end{equation}
where to save space we write
$$
s_x = \sin(\pi x/N),
$$
$$
s_k = \sin(\pi k/N),
$$
$$
f_k = f\big((N/\pi)^2 \sin^2 (\pi k /N) \big),
$$
and
$$
f'_k = f'\big((N/\pi)^2 \sin^2 (\pi k /N) \big).
$$

We always have $0 < k < N/2$, and for fixed $x$ we have $0 < x <
N/2$ whenever $N$ is sufficiently large relative to $x$.  It
follows that $s_k$ and $s_x$ are always within fixed constant
factors of $\pi k/N$ and $\pi x/N$, respectively.  If $k$ is also
sufficiently large relative to $x$, then examining
\eqref{eq:bigformula} and applying $N^2 s_k^2 = O\big(k^2\big)$
shows that the right side is $O\big(k^{-1-2\varepsilon}\big)$ in
absolute value, where the implicit constant is independent of $N$.
\end{proof}

The desired inequalities for $h$ and $\widehat{h}$ now follow
easily.  For each $N$,
$$
h_{N,g}\big(\cos (2\pi x/N)\big) \le g\big(2 - 2\cos (2\pi
x/N)\big),
$$
so it follows in the limit as $N \to \infty$ that
$$
h(x) \le \ft(x),
$$
as desired.  For the other inequality, we apply Bochner's theorem,
which says that a continuous function $h$ on $\R$ is the Fourier
transform of a nonnegative, $L^1$ function (which is of course
uniquely determined by $h$, up to a set of measure $0$) if and
only if $h$ is positive definite, in the sense that for all
$x_1,\dots,x_\ell \in \R$, the $\ell \times \ell$ matrix whose
$i,j$ entry is $h(x_i-x_j)$ is positive semidefinite. (See
\cite[p.~124]{St} or \cite[p.~303]{R}.)  That property follows
immediately from the corresponding property of $h_{N,g}$, because
the set of positive-semidefinite matrices is closed.  This
completes the proof of Proposition~\ref{prop:onedim}.

\section*{Acknowledgements}

We thank Eiichi Bannai, Charles Doran, Noam Elkies, Gabriele Nebe,
Andreas Str\"ombergsson, Frank Vallentin, and G\"unter Ziegler for
their valuable comments and suggestions on our manuscript.

\appendix
\section{Uniqueness of spherical codes}
\label{app}

Uniqueness of the $N$-gon is trivial.  For the simplex, cross
polytope, icosahedron, and $600$-cell, uniqueness follows from
analyzing the case of equality in the B\H or\H oczky bound (see
\cite[p.~260]{B}, or see \cite{BD} for another proof for the
$600$-cell). For the $(24,196560,1/2)$, $(23,4600,1/3)$,
$(8,240,1/2)$, and $(7,56,1/3)$ spherical codes, uniqueness was
proved by Bannai and Sloane in \cite{BS} (see also Chapter~14 of
\cite{CS}).  For the $(6,27,1/4)$, $(5,16,1/5)$, $(22,275,1/6)$,
$(21,162,1/7)$, and $(22,100,1/11)$ codes, uniqueness follows from
results in the theory of strongly regular graphs.  In each case
there are only two distinct inner products other than $1$, so one
can form a graph out of the vertices by making edges correspond to
one particular inner product (one obtains one of two complementary
graphs, which encode the same information).  Using Theorem~7.4 and
Lemma~7.3 from \cite{DGS} one proves that these graphs are
strongly regular and determines their parameters. In each case
there are unique strongly regular graphs with these parameters.
For the $(6,27,1/4)$ and $(5,16,1/5)$ codes, see Theorem~10.6.4
and Lemma~10.9.4 in \cite{GR};  the remaining cases are dealt with
in \cite{CGS}.  It follows that the Gram matrices of the
corresponding point configurations are uniquely determined, so the
configurations are determined up to orthogonal transformations.

For the $(23,552,1/5)$ code, uniqueness follows from the
uniqueness of the regular two-graph on $276$ vertices.  In any
$(23,552,1/5)$ code, the linear programming bounds show that only
the inner products $-1$ and $\pm 1/5$ occur. Such a code thus
gives rise to an arrangement of $276$ equiangular lines in
$\R^{23}$, which corresponds to a regular two-graph on $276$
vertices. (See Chapter~11 of \cite{GR} for background on
two-graphs and equiangular line arrangements.)  Goethals and
Seidel proved in \cite{GS} that it is unique, which implies the
uniqueness of the $(23,552,1/5)$ code.

The last remaining case is the $(22,891,1/4)$ spherical code.  A
proof of its uniqueness is implicit in \cite{Cu} (based on results
from \cite{SY}).  For an alternative proof using the techniques of
Bannai and Sloane, see \cite{CK3}.

\end{document}